\documentclass[a4paper]{amsart}
\usepackage{amssymb}
\usepackage[T1]{fontenc}
\usepackage[latin1]{inputenc}
\usepackage{amsfonts}
\usepackage{amsxtra}
\usepackage{ae}
\usepackage[all]{xy}
\usepackage{enumerate}
 \usepackage[dvips]{color}
\usepackage{color}
\include{diagram}

\newcommand*{\ket}{\rangle}
\newcommand*{\bra}{\langle}
\newcommand*{\ad}{\mathsf{ad}}

\newcommand*{\A}{\mathcal{A}}
\newcommand*{\C}{\mathcal{C}}

\newcommand*{\M}{\mathcal{M}}
\renewcommand*{\M}{\mathcal{M}}
\newcommand*{\N}{\mathcal{N}}
\newcommand*{\E}{\mathcal{E}}
\newcommand*{\F}{\mathcal{F}}

\newcommand*{\X}{\mathcal{X}}
\renewcommand*{\H}{\mathcal{H}}

\newcommand*{\I}{\mathcal{I}}

\newcommand*{\cotimes}{\hat{\otimes}}
\newcommand*{\twisted}{\boxtimes}

\newcommand*{\yd}{\mathsf{YD}}
\newcommand*{\DD}{\mathsf{D}}

\renewcommand*{\max}{\mathsf{f}}
\newcommand*{\red}{\mathsf{r}}

\newcommand*{\cop}{\mathsf{cop}}

\newcommand*{\Alg}{$-$\mathsf{Alg}}
\newcommand*{\hit}{\rightharpoonup}
\newcommand*{\hitby}{\leftharpoonup}

\newcommand*{\HH}{\mathbb{H}}
\newcommand*{\KH}{\mathbb{K}}
\newcommand*{\LH}{\mathbb{L}}

\DeclareMathOperator{\Hom}{Hom}

\DeclareMathOperator{\res}{res}
\DeclareMathOperator{\ind}{ind}

\DeclareMathOperator{\id}{id}
\DeclareMathOperator{\ev}{ev}

\newenvironment{bnum}
{\begin{list}{}
    {\setlength{\labelwidth}{15pt}
     \setlength{\leftmargin}{\labelwidth}
    }
}
{\end{list}}

\numberwithin{equation}{section}
\theoremstyle{change}
\newtheorem{theorem}{Theorem}[section]
\newtheorem{prop}[theorem]{Proposition}
\newtheorem{lemma}[theorem]{Lemma}
\newtheorem{cor}[theorem]{Corollary}
\newtheorem{definition}[theorem]{Definition}

\begin{document}

\title[Poincar\'e duality]{Equivariant Poincar\'e duality for quantum group actions}
\author{Ryszard Nest and Christian Voigt}
\address{Ryszard Nest \\
         Institut for Matematiske Fag \\
         Universitet K\o benhavn \\
         Universitetsparken 5 \\
         2100 K\o benhavn \\
         Denmark
}
\email{rnest@math.ku.dk}
\address{Christian Voigt \\
         Mathematisches Institut\\
         Georg-August-Universit\"at G\"ottingen\\
         Bunsenstra\ss e 3-5 \\
         37073 G\"ottingen\\
         Germany
}
\email{cvoigt@uni-math.gwdg.de}


\maketitle

\begin{abstract}
We extend the notion of Poincar\'e duality in $ KK $-theory to the setting of quantum group actions. An important ingredient in our approach is the replacement of ordinary tensor products by braided tensor products. Along the way we discuss general properties of equivariant $ KK $-theory for locally compact quantum groups, including the construction of exterior products. As an example, we prove that the standard Podle\'s sphere is equivariantly Poincar\'e dual to itself.
\end{abstract}

\section{Introduction}

The notion of Poincar\'e duality in $ K $-theory plays an important r\^ole in noncommutative geometry. 
In particular, it is a fundamental ingredient in the theory of noncommutative manifolds due to Connes \cite{Connes2}. \\
A noncommutative manifold is given by a spectral triple $ (\A, H, D) $ where
$ \A $ is a $ * $-algebra represented on a Hilbert space $ H $ and $ D $ is an unbounded self-adjoint operator on $ H $.
The basic requirements on this data are that $ D $ has compact resolvent and that the commutators $ [D,a] $ are bounded for all $ a \in \A $.
There are further ingredients in the definition of a noncommutative manifold,
in particular a grading and the concept of a real structure \cite{Connesreal}, \cite{Connesgravmat}.
An important recent result due to Connes is the reconstruction theorem
\cite{Connesreconst}, which asserts that in the commutative case, under some natural conditions, the algebra $ \A $ is isomorphic
to $ C^\infty(M) $ for a unique smooth manifold $ M $. The real structure produces a version of $ KO $-Poincar\'e duality, which is a necessary 
ingredient for the existence of a smooth structure. \\
Quantum groups and their homogeneous spaces give natural and interesting examples of noncommutative spaces, and
several cases of associated spectral triples have been constructed \cite{DSPodles}, \cite{CP}, \cite{DLSS}, \cite{DDLW}, \cite{NTDirac}.
An important guiding principle in all these constructions is equivariance with respect to the action of a quantum group.
In \cite{Sitarz}, \cite{Dabrowskispheres} a general framework for equivariant spectral triples is formulated,
including an equivariance condition for real structures. However, in some examples the
original axioms in \cite{Connesreal} are only satisfied up to infinitesimals in this setup \cite{DLSS}, \cite{DDLW}. 
The $ K $-theoretic interpretation of a real structure up to infinitesimals is
unclear.\\
In this paper we introduce a notion of $ K $-theoretic Poincar\'e duality which is particularly
adapted to the symmetry of quantum group actions. More precisely, we generalize the definition of Poincar\'e duality in $ KK $-theory given by Connes \cite{Connes2} to $ C^* $-algebras with a coaction of a quantum group using braided tensor products. Braided tensor products are well-known in the
algebraic approach to quantum groups \cite{Majid}, in our context they are constructed using coactions of
the Drinfeld double of a locally compact quantum group. \\
The example we study in detail is the standard Podle\'s sphere, and we prove that it is equivariantly Poincar\'e dual to itself with respect to
the natural action of $ SU_q(2) $. The Drinfeld double of $ SU_q(2) $, appearing as the symmetry group in this case, is the quantum Lorentz group \cite{PW1}, 
a noncompact quantum group built up out of a compact and a discrete part. We remark that the additional symmetry of the Podle\'s sphere which
is encoded in the discrete part of the quantum Lorentz group is not visible classically. \\
The spectral triple corresponding to the Dirac operator on the standard Podle\'s sphere \cite{DSPodles} can be equipped with
a real structure, and, due to \cite{Wagner}, it satisfies Poincar\'e duality in the sense of
\cite{Connesreal}. From this point of view the standard Podle\'s sphere is very well-behaved. However, already in this example the formulation of
equivariant Poincar\'e duality requires the setup proposed in this paper.\\
Usually, the symmetry of an equivariant spectral triple is implemented by the action of a quantized
universal enveloping algebra. In our approach we have to work with coactions of the quantized algebra of functions instead. Both
descriptions are essentially equivalent, but an advantage of coactions is that the
correct definition of equivariant $ K $-theory and $ K $-homology in this setting is already contained in \cite{BSKK}.
In particular, we do not need to consider constructions of equivariant  $ K $-theory as in \cite{NT}, \cite{Wagner} which do not 
extend to general quantum groups. \\
Let us now describe how the paper is organized. 
In the first part of the paper we discuss some results related to locally compact quantum groups and $ KK $-theory. 
Section \ref{secqg} contains an introduction to locally compact quantum groups, their coactions and associated crossed products.
In particular, we review parts of the foundational work of Vaes on induced coactions \cite{Vaesimprimitivity} which are relevant 
to this paper. In section \ref{secyd} we introduce Yetter-Drinfeld-$ C^* $-algebras and braided tensor products and discuss
their basic properties, including compatiblity with induction and restriction. 
Then, in section \ref{seckkg}, we review the
definition of equivariant $ KK $-theory for quantum groups following
Baaj and Skandalis \cite{BSKK}. In particular, we show that $ KK^G $ for a regular locally compact quantum group $ G $
satisfies a universal property as in the group case.
A new feature in the quantum setting is the construction of exterior products for $ KK^G $. 
The non-triviality of it is related to the fact that a tensor product of two algebras with a coaction of a quantum group does not inherit a natural coaction 
in general, in distinction to the case of a group action. We deal with this problem using braided tensor products. \\
Basic facts concerning $ SU_q(2) $ and the standard Podle\'s sphere $ SU_q(2)/T $ are reviewed in section \ref{secsuq2}.
The main definition and results are contained in section \ref{secfredpodles}, where we introduce the concept of equivariant Poincar\'e duality
with respect to quantum group actions and show that $ SU_q(2)/T $ is equivariantly Poincar\'e dual to itself.
As an immediate consequence we determine the equivariant $ K $-homology of the Podle\'s sphere. \\
Let us make some remarks on notation. We write $ \LH(\E,\F) $ for the space of adjointable
operators between Hilbert $ A $-modules $ \E $ and $ \F $. Moreover $ \KH(\E,\F) $ denotes the
space of compact operators. If $ \E = \F $ we write simply $ \LH(\E) $ and $ \KH(\E) $, respectively.
The closed linear span of a subset $ X $ of a Banach space is denoted by $ [X] $. Depending on the context, the symbol
$ \otimes $ denotes either the tensor product of Hilbert spaces, the minimal tensor product of $ C^* $-algebras, or the tensor product
of von Neumann algebras. For operators on multiple tensor products we use the leg numbering notation. \\
It is a pleasure to thank Uli Kr\"ahmer for interesting discussions on the subject of this paper.
The second author is indebted to Stefaan Vaes for helpful explanations concerning induced coactions and braided tensor products.
A part of this work was done during stays of the authors in Warsaw supported by EU-grant MKTD-CT-2004-509794.
We are grateful to Piotr Hajac for his kind hospitality.

\section{Locally compact quantum groups and their coactions} \label{secqg}

In this section we recall basic definitions and results from the theory of locally compact quantum groups
and fix our notation. For more detailed information we refer to the literature \cite{KVLCQG}, \cite{KVvN}, \cite{Vaesimprimitivity}. \\
Let $ \phi $ be a normal, semifinite and faithful weight on a von Neumann algebra $ M $. We use the standard notation
$$
\M^+_\phi = \{ x \in M_+| \phi(x) < \infty \}, \qquad \N_\phi = \{ x \in M| \phi(x^* x) < \infty \}
$$
and write $ M_*^+ $ for the space of positive normal linear functionals on $ M $. Assume that
$ \Delta: M \rightarrow M \otimes M $ is a normal unital $ * $-homomorphism. The weight $ \phi $ is called left invariant
with respect to $ \Delta $ if
$$
\phi((\omega \otimes \id)\Delta(x)) = \phi(x) \omega(1)
$$
for all $ x \in \M_\phi^+ $ and $ \omega \in M_*^+ $. Similarly one defines the notion of a right invariant
weight.
\begin{definition} \label{defqg}
A locally compact quantum group $ G $ is given by a von Neumann algebra $ L^\infty(G) $ together with a normal unital $ * $-homomorphism
$ \Delta: L^\infty(G) \rightarrow L^\infty(G) \otimes L^\infty(G) $ satisfying the coassociativity relation
$$
(\Delta \otimes \id)\Delta = (\id \otimes \Delta)\Delta
$$
and normal semifinite faithful weights $ \phi $ and $ \psi $ on $ L^\infty(G) $ which are left and right invariant, respectively.
\end{definition}
Our notation for locally compact quantum groups is intended to make clear how ordinary locally compact groups can be viewed as quantum groups.
Indeed, if $ G $ is a locally compact group, then the algebra $ L^\infty(G) $ of essentially bounded measurable functions on $ G $ together with the comultiplication
$ \Delta: L^\infty(G) \rightarrow L^\infty(G) \otimes L^\infty(G) $ given by
$$
\Delta(f)(s,t) = f(st)
$$
defines a locally compact quantum group. The weights $ \phi $ and $ \psi $ are given in this case by left and right Haar measures, respectively. \\
Of course, for a general locally compact quantum group $ G $ the notation $ L^\infty(G) $ is purely formal.
Similar remarks apply to the $ C^* $-algebras $ C^*_\max(G), C^*_\red(G) $ and $ C^\max_0(G), C^\red_0(G) $ associated to $ G $ that we discuss below.
It is convenient to view all of them as different appearances of the quantum group $ G $. \\
Let $ G $ be a locally compact quantum group and let $ \Lambda: \N_\phi \rightarrow \HH_G $ be a GNS-construction for the weight $ \phi $.
Throughout the paper we will only consider quantum groups for which $ \HH_G $ is a separable Hilbert space. One obtains a unitary
$ W_G = W $ on $ \HH_G \otimes \HH_G $ by
$$
W^*(\Lambda(x) \otimes \Lambda(y)) = (\Lambda \otimes \Lambda)(\Delta(y)(x \otimes 1))
$$
for all $ x, y \in \N_\phi $. This unitary is multiplicative, which means that $ W $ satisfies the pentagonal equation
$$
W_{12} W_{13} W_{23} = W_{23} W_{12}.
$$
From $ W $ one can recover the von Neumann algebra $ L^\infty(G) $ as the strong closure of the algebra
$ (\id \otimes \LH(\HH_G)_*)(W) $ where $ \LH(\HH_G)_* $ denotes the space of normal linear functionals on $ \LH(\HH_G) $. Moreover
one has
$$
\Delta(x) = W^*(1 \otimes x) W
$$
for all $ x \in M $. The algebra $ L^\infty(G) $ has an antipode which
is an unbounded, $ \sigma $-strong* closed linear map $ S $ given by $ S(\id \otimes \omega)(W) = (\id \otimes \omega)(W^*) $
for $ \omega \in \LH(\HH_G)_* $. Moreover there is a polar decomposition $ S = R \tau_{-i/2} $ where $ R $ is an
antiautomorphism of $ L^\infty(G) $ called the unitary antipode and $ (\tau_t) $ is a strongly continuous one-parameter group
of automorphisms of $ L^\infty(G) $ called the scaling group. The unitary antipode satisfies $ \sigma(R \otimes R) \Delta = \Delta R $. \\
The group-von Neumann algebra $ \mathcal{L}(G) $ of the quantum group $ G $ is the strong
closure of the algebra $ (\LH(\HH_G)_* \otimes \id)(W) $ with the comultiplication $ \hat{\Delta}: \mathcal{L}(G) \rightarrow \mathcal{L}(G) \otimes \mathcal{L}(G) $
given by
$$
\hat{\Delta}(y) = \hat{W}^* (1 \otimes y) \hat{W}
$$
where $ \hat{W} = \Sigma W^* \Sigma $ and $ \Sigma \in \LH(\HH_G \otimes \HH_G) $ is the flip map. It defines
a locally compact quantum group $ \hat{G} $ which is called the dual of $ G $. The left invariant weight
$ \hat{\phi} $ for the dual quantum group has a GNS-construction $ \hat{\Lambda}: \N_{\hat{\phi}} \rightarrow \HH_G $,
and according to our conventions we have $ \mathcal{L}(G) = L^\infty(\hat{G}) $. \\
The modular conjugations of the weights $ \phi $ and $ \hat{\phi} $ are denoted by $ J $ and $ \hat{J} $, respectively.
These operators implement the unitary antipodes in the sense that
$$
R(x) = \hat{J} x^* \hat{J}, \qquad \hat{R}(y) = J y^* J
$$
for $ x \in L^\infty(G) $ and $ y \in \mathcal{L}(G) $. Note that $ L^\infty(G)' = JL^\infty(G) J $ and
$ \mathcal{L}(G)' = \hat{J} \mathcal{L}(G) \hat{J} $ for the commutants of $ L^\infty(G) $ and $ \mathcal{L}(G) $. Using $ J $ and $ \hat{J} $ one obtains multiplicative
unitaries
$$
V = (\hat{J} \otimes \hat{J})\hat{W}(\hat{J} \otimes \hat{J}), \qquad \hat{V} = (J \otimes J) W (J \otimes J)
$$
which satisfy $ V \in \mathcal{L}(G)' \otimes L^\infty(G) $ and $ \hat{V} \in L^\infty(G)' \otimes \mathcal{L}(G) $, respectively. \\
We will mainly work with the $ C^* $-algebras associated to the locally compact quantum group $ G $. The
algebra $ [(\id \otimes \LH(\HH_G)_*)(W)] $ is a strongly dense $ C^* $-subalgebra of $ L^\infty(G) $
which we denote by $ C^\red_0(G) $. Dually, the algebra
$ [(\LH(\HH_G)_* \otimes \id)(W)] $ is a strongly dense $ C^* $-subalgebra of $ \mathcal{L}(G) $
which we denote by $ C^*_\red(G) $.
These algebras are the reduced algebra of continuous functions vanishing at infinity
on $ G $ and the reduced group $ C^* $-algebra of $ G $, respectively. One has
$ W \in M(C^\red_0(G) \otimes C^*_\red(G)) $. \\
Restriction of the comultiplications on $ L^\infty(G) $ and $ \mathcal{L}(G) $
turns $ C^\red_0(G) $ and $ C^*_\red(G) $ into Hopf-$ C^* $-algebras in the following sense.
\begin{definition}\label{defhopfcstar}
A Hopf $ C^* $-algebra is a $ C^* $-algebra $ S $ together with an injective nondegenerate $ * $-homomorphism
$ \Delta: S \rightarrow M(S \otimes S) $ such that the diagram
$$
\xymatrix{
S \ar@{->}[r]^{\Delta} \ar@{->}[d] & M(S \otimes S) \ar@{->}[d]^{\id \otimes \Delta} \\
M(S \otimes S) \ar@{->}[r]^{\!\!\!\!\!\!\!\!\! \Delta \otimes \id} & M(S \otimes S \otimes S)
     }
$$
is commutative and $ [\Delta(S)(1 \otimes S)] = S \otimes S = [(S \otimes 1)\Delta(S)] $. \\
A morphism between Hopf-$ C^* $-algebras $ (S, \Delta_S) $ and $ (T, \Delta_T) $ is a nondegenerate $ * $-homomorphism
$ \pi: S \rightarrow M(T) $ such that $ \Delta_T\, \pi = (\pi \otimes \pi)\Delta_S $.
\end{definition}
If $ S $ is a Hopf-$ C^* $-algebra we write $ S^\cop $ for the Hopf-$ C^* $-algebra obtained by equipping
$ S $ with the opposite comultiplication $ \Delta^\cop = \sigma \Delta $. \\
A unitary corepresentation of a Hopf-$ C^* $-algebra $ S $ on a Hilbert $ B $-module $ \E $ is
a unitary $ X \in \LH(S \otimes \E) $ satisfying
$$
(\Delta \otimes \id)(X) = X_{13} X_{23}.
$$
A universal dual of $ S $ is a Hopf-$ C^* $-algebra $ \hat{S} $ together with a unitary
corepresentation $ \X \in M(S \otimes \hat{S}) $ satisfying the following universal property.
For every Hilbert $ B $-module $ \E $
and every unitary corepresentations $ X \in \LH(S \otimes \E) $ there exists a unique nondegenerate $ * $-homomorphism
$ \pi_X: \hat{S} \rightarrow \LH(\E) $ such that $ (\id \otimes \pi_X)(\X) = X $. \\
For every locally compact quantum group $ G $ there exists a universal dual $ C^*_\max(G) $ of $ C_0^\red(G) $ and
a universal dual $ C^\max_0(G) $ of $ C^*_\red(G) $, respectively \cite{Kustermansuniversal}.
We call $ C^*_\max(G) $ the maximal group $ C^* $-algebra of $ G $ and
$ C_0^\max(G) $ the maximal algebra of continuous functions on $ G $ vanishing at infinity.
Since $ \HH_G $ is assumed to be separable the $ C^* $-algebras $ C^\max_0(G), C^\red_0(G) $ and
$ C^*_\max(G), C^*_\red(G) $ are separable. The quantum group $ G $ is called compact if $ C^\max_0(G) $ is unital, and it is called
discrete if $ C^*_\max(G) $ is unital. In the compact case we also write $ C^\max(G) $ and $ C^\red(G) $ instead of $ C^\max_0(G) $ and $ C^\red_0(G) $,
respectively. \\
In general, we have a surjective morphism $ \hat{\pi}: C^*_\max(G) \rightarrow C^*_\red(G) $ of Hopf-$ C^* $-algebras
associated to the left regular corepresentation $ W \in M(C_0(G) \otimes C^*_\red(G)) $. Similarly, there is
a surjective morphism $ \pi: C^\max_0(G) \rightarrow C^\red_0(G) $.
We will call the quantum group $ G $ amenable if $ \hat{\pi}: C^*_\max(G) \rightarrow C^*_\red(G) $
is an isomorphism and coamenable if $ \pi: C^\max_0(G) \rightarrow C^\red_0(G) $
is an isomorphism. If $ G $ is amenable or coamenable, respectively, we also write $ C^*(G) $ and $ C_0(G) $
for the corresponding $ C^* $-algebras. For more information on amenability for locally compact quantum groups see \cite{BedosTuset}. \\
Let $ S $ be a $ C^* $-algebra. The $ S $-relative multiplier algebra $ M_S(S \otimes A) \subset M(S \otimes A) $ of a $ C^* $-algebra $ A $
consists of all $ x \in M(S \otimes A) $ such that the relations
$$
x(S \otimes 1) \subset S \otimes A, \qquad (S \otimes 1)x \subset S \otimes A
$$
hold. In the sequel we tacitly use basic properties of relative multiplier algebras which can be found in \cite{EKQR}.
\begin{definition} \label{defcoaction}
A (left) coaction of a Hopf $ C^* $-algebra $ S $ on a $ C^* $-algebra $ A $ is an injective nondegenerate $ * $-homomorphism
$ \alpha: A \rightarrow M(S \otimes A) $ such that the diagram
$$
\xymatrix{
A \ar@{->}[r]^{\alpha} \ar@{->}[d]^\alpha & M(S \otimes A) \ar@{->}[d]^{\Delta \otimes \id} \\
M(S \otimes A) \ar@{->}[r]^{\!\!\!\!\!\!\!\!\!\! \id \otimes \alpha} & M(S \otimes S \otimes A)
     }
$$
is commutative and $ \alpha(A) \subset M_S(S \otimes A) $.
The coaction is called continuous if $ [\alpha(A)(S \otimes 1)] = S \otimes A $. \\
If $ (A, \alpha) $ and $ (B, \beta) $ are $ C^* $-algebras with coactions of $ S $, then a $ * $-homomorphism
$ f: A \rightarrow M(B) $ is called $ S $-colinear if $ \beta f = (\id \otimes f)\alpha $.
\end{definition}
We remark that some authors do not require a coaction to be injective. For a discussion of the
continuity condition see \cite{BSV}. \\
A $ C^* $-algebra $ A $ equipped with a continuous coaction of the Hopf-$ C^* $-algebra $ S $ will be called an $ S $-$ C^* $-algebra.
If $ S = C^\red_0(G) $ for a locally compact quantum group $ G $ we also say that $ A $ is
$ G $-$ C^* $-algebra. Moreover, in this case $ S $-colinear $ * $-homomorphisms will be called $ G $-equivariant or simply equivariant.
We write $ G \Alg $ for the category of separable $ G $-$ C^* $-algebras and
equivariant $ * $-homomorphisms. \\
A (nondegenerate) covariant representation of a $ G $-$ C^* $-algebra $ A $ on a Hilbert-$ B $-module $ \E $
consists of a (nondegenerate) $ * $-homomorphism $ f: A \rightarrow \LH(\E) $ and a unitary
corepresentation $ X \in \LH(C^\red_0(G) \otimes \E) $ such that
$$
(\id \otimes f)\alpha(a) = X^*(1 \otimes f(a)) X
$$
for all $ a \in A $. There exists a $ C^* $-algebra $ C^*_\max(G)^\cop \ltimes_\max A $, called the full crossed product,
together with a nondegenerate covariant representation $ (j_A, X_A) $ of $ A $ on $ C^*_\max(G)^\cop \ltimes_\max A $
which satisfies the following universal property.
For every nondegenerate covariant representation $ (f, X) $ of $ A $ on a Hilbert-$ B $-module $ \E $ there exists a unique
nondegenerate $ * $-homomorphism $ F: C^*_\max(G)^\cop \ltimes_\max A \rightarrow \LH(\E) $, called the integrated form
of $ (f,X) $, such that
$$
X = (\id \otimes F)(X_A), \qquad f = F j_A.
$$
Remark that the corepresentation $ X_A $ corresponds to a unique nondegenerate $ * $-homomorphism
$ g_A: C^*_\max(G)^\cop \rightarrow M(C^*_\max(G)^\cop \ltimes_\max A) $. \\
On the Hilbert $ A $-module $ \HH_G \otimes A $ we have a covariant representation of $ A $
given by the coaction $ \alpha: A \rightarrow \LH(\HH_G \otimes A) $ and $ W \otimes 1 $.
The reduced crossed product $ C^*_\red(G)^\cop \ltimes_\red A $ is the image of $ C^*_\max(G)^\cop \ltimes_\max A $
under the corresponding integrated form. Explicitly, we have
$$
C^*_\red(G)^\cop \ltimes_\red A = [(C^*_\red(G) \otimes 1)\alpha(A)]
$$
inside $ M(\KH_G \otimes A) = \LH(\HH_G \otimes A) $ using the notation $ \KH_G = \KH(\HH_G) $.
There is a nondegenerate $ * $-homomorphism $ j_A: A \rightarrow M(C^*_\red(G)^\cop \ltimes_\red A) $ induced by $ \alpha $.
Similarly, we have a canonical nondegenerate $ * $-homomorphism $ g_A: C^*_\red(G)^\cop \rightarrow M(C^*_\red(G)^\cop \ltimes_\red A) $. \\
The full and the reduced crossed products admit continuous dual coactions of $ C^*_\max(G)^\cop $ and $ C^*_\red(G)^\cop $, respectively.
In both cases the dual coaction leaves the copy of $ A $ inside the crossed product invariant and acts by the (opposite) comultiplication on
the group $ C^* $-algebra. If $ G $ is amenable then the canonical map $ C^*_\max(G)^\cop \ltimes_\max A \rightarrow C^*_\red(G)^\cop \ltimes_\red A $ is
an isomorphism for all $ G $-$ C^* $-algebras $ A $, and we will also write $ C^*(G)^\cop \ltimes A $ for the crossed product in this case. \\
The comultiplication $ \Delta: C^\red_0(G) \rightarrow M(C^\red_0(G) \otimes C^\red_0(G)) $ defines a coaction of $ C^\red_0(G) $
on itself. On the Hilbert space $ \HH_G $ we have a covariant representation of $ C^\red_0(G) $ given by the
identical representation of $ C^\red_0(G) $ and $ W \in M(C^\red_0(G) \otimes \KH_G) $.
The quantum group $ G $ is called strongly regular if the associated integrated form induces an isomorphism
$ C^*_\max(G)^\cop \ltimes_\max C^\red_0(G) \cong \KH_G $. Similarly, $ G $ is called regular if the corresponding homomorphism
on the reduced level gives an isomorphism $ C^*_\red(G)^\cop \ltimes_\red C^\red_0(G) \cong \KH_G $.
Every strongly regular quantum group is regular, it is not known wether there exist regular quantum groups which are not strongly regular. If $ G $ is regular then the dual $ \hat{G} $ is regular as well. \\
Let $ \E_B $ be a right Hilbert module. The multiplier module $ M(\E) $ of $ \E $ is the right Hilbert-$ M(B) $-module
$ M(\E) = \LH(B, \E) $. There is a natural embedding $ \E \cong \KH(B, \E) \rightarrow \LH(B, \E) = M(\E) $. If $ \E_A $ and $ \F_B $ are
Hilbert modules, then a morphism from $ \E $ to $ \F $ is a linear map $ \Phi: \E \rightarrow M(\F) $ together with a $ * $-homomorphism
$ \phi: A \rightarrow M(B) $ such that
$$
\bra \Phi(\xi), \Phi(\eta) \ket = \phi(\bra \xi, \eta \ket)
$$
for all $ \xi, \eta \in \E $. In this case $ \Phi $ is automatically norm-decreasing and satisfies
$ \Phi(\xi a) = \Phi(\xi) \phi(a) $ for all $ \xi \in \E $ and $ a \in A $.
The morphism $ \Phi $ is called nondegenerate if $ \phi $ is nondegenerate and $ [\Phi(\E) B] = \F $. \\
Let $ S $ be a $ C^* $-algebra and let $ \E_A $ be a Hilbert module.
The $ S $-relative multiplier module $ M_S(S \otimes \E) $ is the submodule of $ M(S \otimes \E) $
consisting of all multipliers $ x $ satisfying $ x(S \otimes 1) \subset S \otimes \E $ and $ (S \otimes 1)x \subset S \otimes \E $.
For further information we refer again to \cite{EKQR}.
\begin{definition}\label{defcoactionhilbert}
Let $ S $ be a Hopf-$ C^* $-algebra and let $ \beta: B \rightarrow M(S \otimes B) $ be a coaction of $ S $ on the $ C^* $-algebra $ B $.
A coaction of $ S $ on a Hilbert module $ \E_B $ is a nondegenerate morphism $ \lambda: \E \rightarrow M(S \otimes \E) $ such that the
diagram
$$
\xymatrix{
\E \ar@{->}[r]^\lambda \ar@{->}[d]^\lambda & M(S \otimes \E) \ar@{->}[d]^{\Delta \otimes \id} \\
M(S \otimes \E) \ar@{->}[r]^{\!\!\!\!\!\!\!\!\!\! \id \otimes \lambda} & M(S \otimes S \otimes \E)
     }
$$
is commutative and $ \lambda(\E) \subset M_S(S \otimes \E) $. The coaction $ \lambda $ is called continuous if
$ [(S \otimes 1)\lambda(\E)] = S \otimes \E = [\lambda(\E)(S \otimes 1)] $. \\
A morphism $ \Phi: \E \rightarrow M(\F) $ of Hilbert $ B $-modules with coactions $ \lambda_\E $ and $ \lambda_\F $, respectively,
is called $ S $-colinear if $ \lambda_\F \Phi = (\id \otimes \Phi)\lambda_\E $.
\end{definition}
If $ \lambda: \E \rightarrow M(S \otimes \E) $ is a coaction on the Hilbert-$ B $-module $ \E $ then the map $ \lambda $ is automatically isometric and
hence injective. \\
Let $ G $ be a locally compact quantum group and let $ B $ be a $ G $-$ C^* $-algebra. A $ G $-Hilbert $ B $-module
is a Hilbert module $ \E_B $ with a continuous coaction $ \lambda: \E \rightarrow M(S \otimes \E) $ for $ S = C^\red_0(G) $.
If $ G $ is regular then continuity of the coaction $ \lambda $ is in fact automatic. Instead of
$ S $-colinear morphisms we also speak of equivariant morphisms between $ G $-Hilbert $ B $-modules. \\
Let $ B $ be a $ C^* $-algebra equipped with a coaction of the Hopf-$ C^* $-algebra $ S $. Given a Hilbert module $ \E_B $ with coaction
$ \lambda: \E \rightarrow M(S \otimes \E) $ one obtains
a unitary operator $ V_\lambda: \E \otimes_B (S \otimes B) \rightarrow S \otimes \E $ by
$$
V_\lambda(\xi \otimes x) = \lambda(\xi) x
$$
for $ \xi \in \E $ and $ x \in S \otimes B $. Here the tensor product over $ B $ is formed with respect to the
coaction $ \beta: B \rightarrow M(S \otimes B) $. This unitary satisfies the
relation
$$
(\id \otimes_\mathbb{C} V_\lambda)(V_\lambda \otimes_{(\id \otimes \beta)} \id) = V_\lambda \otimes_{(\Delta \otimes \id)} \id
$$
in $ \LH(\E \otimes_{(\Delta \otimes \id)\beta} (S \otimes S \otimes B), S \otimes S \otimes \E) $, compare \cite{BSKK}.
Moreover, the equation
$$
\ad_\lambda(T) = V_\lambda (T \otimes \id) V_\lambda^*
$$
determines a coaction $ \ad_\lambda: \KH(\E) \rightarrow M(S \otimes \KH(\E)) = \LH(S \otimes \E) $.
If the coaction $ \lambda $ is continuous then $ \ad_\lambda $ is continuous as well.
In particular, if $ \E $ is a $ G $-Hilbert $ B $-module with coaction $ \lambda $,
then the associated coaction $ \ad_\lambda $ turns $ \KH(\E) $ into a $ G $-$ C^* $-algebra. \\
Let $ B $ be a $ C^* $-algebra equipped with the trivial coaction of the Hopf-$ C^* $-algebra $ S $ and let $ \lambda: \E \rightarrow M(S \otimes \E) $
be a coaction on the Hilbert module $ \E_B $. Then using the natural identification
$ \E \otimes_B (S \otimes B) \cong \E \otimes S \cong S \otimes \E $ the associated unitary $ V_\lambda $ determines
a unitary corepresentation $ V_\lambda^* $ in $ \LH(S \otimes \E)  $.
Conversely, if $ V \in \LH(S \otimes \E) $ is a unitary corepresentation then
$ \lambda_V: \E \rightarrow M(S \otimes \E) $ given by $ \lambda_V(\xi) = V^*(1 \otimes \xi) $ is
a nondegenerate morphism of Hilbert modules satisfying the coaction identity.
If $ S = C^\red_0(G) $ for a regular quantum group $ G $, then $ \lambda_V $ defines a continuous
coaction on $ \E $. As a consequence, for a regular quantum group $ G $ and a trivial $ G $-$ C^* $-algebra $ B $,
continuous coactions on a Hilbert $ B $-module $ \E $ correspond uniquely to
unitary corepresentations of $ C^\red_0(G) $ on $ \E $. \\
Let $ G $ be a regular quantum group and let $ \E_B $ be a $ G $-Hilbert module with coaction $ \lambda_\E: \E \rightarrow M(C^\red_0(G) \otimes \E) $.
Then $ \HH_G \otimes \E $ becomes a $ G $-Hilbert $ B $-module with the coaction $ \lambda_{\HH_G \otimes \E}(x \otimes \xi)
= X_{12}^* \Sigma_{12} (\id \otimes \lambda_\E)(x \otimes \xi) $ where $ X = \Sigma V \Sigma \in \LH(C^\red_0(G) \otimes \KH_G) $.
In particular, for $ \E = B $ the algebra $ \KH_G \otimes B = \KH(\HH_G \otimes B) $ can be viewed as a $ G $-$ C^* $-algebra.
We now state the following version of the Takesaki-Takai duality theorem \cite{BSUM}.
\begin{theorem} \label{TTduality}
Let $ G $ be a regular locally compact quantum group and let $ A $ be a $ G $-$ C^* $-algebra. Then
there is a natural isomorphism
$$
C^\red_0(G) \ltimes_\red C^*_\red(G)^\cop \ltimes_\red A \cong \KH_G \otimes A
$$
of $ G $-$ C^* $-algebras.
\end{theorem}
An equivariant Morita equivalence between $ G $-$ C^* $-algebras $ A $ and $ B $ is given by an equivariant $ A $-$ B $-imprimitivity bimodule, that is,
a full $ G $-Hilbert $ B $-module $ \E $ together with an isomorphism
$ A \cong \KH(\E) $ of $ G $-$ C^* $-algebras. Theorem \ref{TTduality} shows that the double crossed product
$ C^\red_0(G) \ltimes_\red C^*_\red(G)^\cop \ltimes_\red A $ is equivariantly Morita equivalent to $ A $ for every $ G $-$ C^* $-algebra $ A $
provided $ G $ is regular. \\
A morphism $ H \rightarrow G $ of locally compact quantum groups is a nondegenerate $ * $-homomorphism $ \pi: C^\max_0(G) \rightarrow M(C^\max_0(H)) $
which is compatible with the comultiplications in the sense that $ (\pi \otimes \pi) \Delta_G = \Delta_H \pi $.
Every such morphism induces canonically a dual morphism $ \hat{\pi}: C^*_\max(H) \rightarrow M(C^*_\max(G)) $.
A closed quantum subgroup $ H \subset G $ is a morphism $ H \rightarrow G $ for which the latter map is accompanied by a faithful normal $ * $-homomorphism
$ \mathcal{L}(H) \rightarrow \mathcal{L}(G) $ of\textsc{} the group-von Neumann algebras, see \cite{VaesVainermanlowdim}, \cite{Vaesimprimitivity}.
In the classical case this notion recovers precisely the closed subgroups of a locally compact group $ G $.
Observe that there is in general no associated homomorphism $ L^\infty(G) \rightarrow L^\infty(H) $ for a quantum subgroup, this fails already in
the group case. \\
Let $ H \rightarrow G $ be a morphism of quantum groups and let $ B $ be a $ G $-$ C^* $-algebra with coaction 
$ \beta: B \rightarrow M(C^\red_0(G) \otimes B) $.
Identifying $ \beta $ with a normal coaction \cite{Fischerthesis} of the full $ C^* $-algebra $ C^\max_0(G) $, the map
$ \pi: C^\max_0(G) \rightarrow M(C^\max_0(H)) $ induces on $ B $ a continuous coaction $ \res(\beta): B \rightarrow M(C^\red_0(H) \otimes B) $.
We write $ \res^G_H(B) $ for the resulting $ H $-$ C^* $-algebra. In this way we obtain a functor $ \res^G_H: G \Alg \rightarrow H \Alg $. \\
Conversely, let $ G $ be a strongly regular quantum group and let $ H \subset G $ be a closed quantum subgroup.
Given an $ H $-$ C^*$-algebra $ B $, there exists an induced $ G $-$ C^* $-algebra
$ \ind_H^G(B) $ such that the following version of Green's imprimitivity theorem holds \cite{Vaesimprimitivity}.
\begin{theorem}\label{vaesimprimitivity}
Let $ G $ be a strongly regular quantum group and let $ H \subset G $ be a closed quantum subgroup. Then there is a natural
$ C^*_\red(G)^\cop $-colinear Morita equivalence
$$
C^*_\red(G)^\cop \ltimes_\red \ind_H^G(B) \sim_M C^*_\red(H)^\cop \ltimes_\red B
$$
for all $ H $-$ C^* $-algebras $ B $.
\end{theorem}
In fact, the induced $ C^* $-algebra $ \ind_H^G(B) $ is defined by Vaes in \cite{Vaesimprimitivity} using a generalized Landstad theorem after construction of its reduced crossed product.
A description of $ C^*_\red(G)^\cop \ltimes_\red \ind_H^G(B) $ can be given as follows. From the quantum subgroup $ H \subset G $
one first obtains a right coaction $ L^\infty(G) \rightarrow L^\infty(G) \otimes L^\infty(H) $ on the level of von Neumann algebras.
The von-Neumann algebraic homogeneous space $ L^\infty(G/H) \subset L^\infty(G) $ is defined as the subalgebra of invariants under this
coaction. If $ \hat{\pi}': \mathcal{L}(H)' \rightarrow \mathcal{L}(G)' $ is the homomorphism $ \hat{\pi}'(x) = \hat{J}_G \hat{\pi}(\hat{J}_H x \hat{J}_H) \hat{J}_G $
induced by $ \hat{\pi}: \mathcal{L}(H) \rightarrow \mathcal{L}(G) $, then
$$
I = \{v \in \LH(\HH_H, \HH_G) |\, v x = \hat{\pi}'(x) v \;\text{for all}\; x \in \mathcal{L}(H)' \}
$$
defines a von-Neumann algebraic imprimitivity bimodule between the von Neumann algebraic crossed product $ \mathcal{L}(G)^\cop \ltimes L^\infty(G/H) $ and $ \mathcal{L}(H)^\cop $.
There is a $ C^* $-algebraic homogeneous space $ C^\red_0(G/H) \subset L^\infty(G/H) $ and
a $ C^* $-algebraic imprimitivity bimodule $ \I \subset I $ which implements a Morita equivalence
between $ C^*_\red(G)^\cop \ltimes_\red C^\red_0(G/H) $ and $ C^*_\red(H)^\cop $.
Explicitly, we have
$$
\I \otimes \HH_G = [\hat{V}_G (I \otimes 1)(\id \otimes \hat{\pi})(\hat{V}_H^*)(C^*_\red(H) \otimes \HH_G)].
$$
The crossed product of the induced $ C^* $-algebra $ \ind_H^G(B) $ is then given by
$$
C^*_\red(G)^\cop \ltimes_\red \ind_H^G(B) = [(\I \otimes 1) \beta(B) (\I^* \otimes 1)]
$$
where $ \beta: B \rightarrow M(C^\red_0(H) \otimes B) $ is the coaction on $ B $. \\
At several points of the paper we will rely on techniques developed in \cite{Vaesimprimitivity}. Firstly, as indicated in \cite{Vaesimprimitivity},
let us note that we have induction in stages.
\begin{prop} \label{inductionstages}
Let $ H \subset K \subset G $ be strongly regular quantum groups. Then there is a natural $ G $-equivariant isomorphism
$$
\ind_H^G(B) \cong \ind_K^G \ind_H^K(B)
$$
for every $ H $-$ C^* $-algebra $ B $.
\end{prop}
\proof Let $ \hat{\pi}_H^G: \mathcal{L}(H) \rightarrow \mathcal{L}(G) $ be the normal $ * $-homomorphism corresponding
to the inclusion $ H \subset G $, and denote by $ \I_H^G \subset I_H^G \subset \LH(\HH_H, \HH_G) $ the associated
imprimitivity bimodules. For the inclusions $ H \subset K $ and $ K \subset G $ we use
analogous notation. By assumption we have $ \hat{\pi}_K^G \hat{\pi}_H^K = \hat{\pi}_H^G $, and
we observe that $ I_K^G I_H^K \subset I_H^G $ is strongly dense. \\
Since the $ * $-homomorphism $ \hat{\pi}_K^G $ is normal and injective we obtain
$$
\I_H^K \otimes \HH_G = [(\id \otimes \hat{\pi}_K^G)(\hat{V}_K)(I_H^K \otimes 1)(\id \otimes \hat{\pi}_H^G)(\hat{V}_H^*)(C^*_\red(H) \otimes \HH_G)]
$$
which yields
$$
[\I_K^G \I_H^K] \otimes \HH_G =
[\hat{V}_G(I_K^G I_H^K \otimes 1)(\id \otimes \hat{\pi}_H^G)(\hat{V}_H^*)(C^*_\red(H) \otimes \HH_G)].
$$
Using the normality of $ \hat{\pi}_H^G $ we see that if $ (v_i)_{i\in I} $ is a bounded net in $ I_H^G $ converging strongly to zero then
$ \hat{V}_G(v_i \otimes 1)(\id \otimes \hat{\pi}_H^G)(\hat{V}_H^*)(x \otimes \xi) $
converges to zero in norm for all $ x \in C^*_\red(H) $ and $ \xi \in \HH_G $.
As a consequence we obtain $ \I_H^G = [\I_K^G \I_H^K] $ for the $ C^* $-algebraic imprimitivity bimodules. \\
Now let $ B $ be an $ H $-$ C^* $-algebra with coaction $ \beta $. Then we have
\begin{align*}
C^*_\red(G)^\cop \ltimes_\red &\ind_K^G \ind_H^K(B) = [(\I_K^G \otimes \id \otimes \id)(\Delta_K \otimes \id)(\ind_H^K(B)) ((\I_K^G)^* \otimes \id \otimes \id)] \\
&\cong [(\I_K^G \otimes \id \otimes \id)(V_K)^*_{12}(\Delta_K \otimes \id)(\ind_H^K(B)) (V_K)_{12}((\I_K^G)^* \otimes \id \otimes \id)] \\
&\cong [(\I_K^G \I_H^K \otimes \id)\beta(B) ((\I_H^K)^* (\I_K^G)^* \otimes \id)] \\
&= [(\I_H^G \otimes \id) \beta(B) ((\I_H^G)^* \otimes \id)] = C^*_\red(G)^\cop \ltimes_\red \ind_H^G(B)
\end{align*}
using conjugation with the unitary $ ((\hat{\pi}_K^{G})' \otimes \id)(V_K^*)_{12} $ in the second step.
The resulting isomorphism between the crossed products $ C^*_\red(G)^\cop \ltimes_\red \ind_K^G \ind_H^K(B) $ and $ C^*_\red(G)^\cop \ltimes_\red \ind_H^G(B) $
is $ C^*_\red(G)^\cop $-colinear and identifies the natural corepresentations of $ C^\red_0(G) $ on both sides.
Hence theorem 6.7 in \cite{Vaesimprimitivity} yields the assertion. \qed \\
Let $ H \subset G $ be a quantum subgroup of a strongly regular quantum group $ G $ and let $ B $ be an $ H $-$ C^* $-algebra
with coaction $ \beta $. If $ E $ denotes the trivial group, then due to proposition \ref{inductionstages} we have
$$
\ind_H^G(C^\red_0(H) \otimes B) = \ind_H^G \ind_E^H \res^H_E(B) \cong \ind_E^G \res^H_E(B) = C^\red_0(G) \otimes B
$$
where $ C^\red_0(H) \otimes B $ is viewed as an $ H $-$ C^* $-algebra via comultiplication on the first tensor factor.
The $ * $-homomorphism $ \beta: B \rightarrow M(C^\red_0(H) \otimes B) $ induces an injective $ G $-equivariant
$ * $-homomorphism $ \ind(\beta): \ind_H^G(B) \rightarrow M(\ind_H^G(C^\red_0(H) \otimes B)) $, and it follows that
$ \ind_H^G(B) $ is contained in $ M(C^\red_0(G) \otimes B) $. Using that the coaction $ \beta $ is continuous
we see that $ \ind_H^G(B) $ is in fact contained in the $ C^\red_0(G) $-relative multiplier algebra of $ C^\red_0(G) \otimes B $. \\
Now let $ A $ and $ B $ be $ H $-$ C^* $-algebras. According to the previous observations every $ H $-equivariant
$ * $-homomorphism $ f: A \rightarrow B $ induces a
$ G $-equivariant $ * $-homomorphism $ \ind_H^G(f): \ind_H^G(A) \rightarrow \ind_H^G(B) $
in a natural way. We conclude that induction defines a functor $ \ind_H^G: H \Alg \rightarrow G \Alg $.

\section{Yetter-Drinfeld algebras and braided tensor products} \label{secyd}

In this section we study Yetter-Drinfeld-$ C^* $-algebras and braided tensor products.
We remark that these concepts are well-known in the algebraic approach to quantum groups \cite{Majid}. Yetter-Drinfeld modules for
compact quantum groups are discussed in \cite{PW1}. \\
Let us begin with the definition of a Yetter-Drinfeld $ C^* $-algebra.
\begin{definition} \label{defyd}
Let $ G $ be a locally compact quantum group and let $ S = C^{\red}_0(G) $ and $ \hat{S} = C^*_\red(G) $ be the
associated reduced Hopf-$ C^* $-algebras. A $ G $-Yetter-Drinfeld $ C^* $-algebra is a
$ C^* $-algebra $ A $ equipped with continuous coactions $ \alpha $ of $ S $ and $ \lambda $ of $ \hat{S} $ such that the diagram
$$
\xymatrix{
A \ar@{->}[r]^\lambda \ar@{->}[d]^\alpha & M(\hat{S} \otimes A) \ar@{->}[r]^{\!\!\!\!\!\!\!\!\! \id \otimes \alpha} &
M(\hat{S} \otimes S \otimes A) \ar@{->}[d]^{\sigma \otimes \id} \\
M(S \otimes A) \ar@{->}[r]^{\!\!\!\!\!\!\!\!\!\! \id \otimes \lambda} & M(S \otimes \hat{S} \otimes A)\; \ar@{->}[r]^{\ad(W) \otimes \id}
& \;M(S \otimes \hat{S} \otimes A)
     }
$$
is commutative. Here $ \ad(W)(x) = WxW^* $ denotes the adjoint action of the fundamental unitary $ W \in M(S \otimes \hat{S}) $.
\end{definition}
In order to compare definition \ref{defyd} with the notion
of a Yetter-Drinfeld module in the algebraic setting, one should keep in mind that we work with the opposite comultiplication
on the dual. In the sequel we will also refer to $ G $-Yetter-Drinfeld $ C^* $-algebras as $ G $-$ \yd $-algebras.
A homomorphism of $ G $-$ \yd $-algebras $ f: A \rightarrow B $ is a $ * $-homomorphism which is both $ G $-equivariant and $ \hat{G} $-equivariant.
We remark moreover that the concept of a Yetter-Drinfeld-$ C^* $-algebra is self-dual, that is, $ G $-$ \yd $-algebras
are the same thing as $ \hat{G} $-$ \yd $-algebras. \\
Let us discuss some basic examples of Yetter-Drinfeld-$ C^* $-algebras. Consider first the case that $ G $ is an ordinary locally compact group.
Since $ C_0(G) $ is commutative, every $ G $-$ C^* $-algebra becomes a $ G $-$ \yd $-algebra with the trivial coaction of $ C^*_\red(G) $.
Dually, we may start with a $ \hat{G} $-$ C^* $-algebra, that is, a reduced coaction of the group $ G $.
If $ G $ is discrete then such coactions correspond to Fell bundles over $ G $. In this case a Yetter-Drinfeld structure
is determined by an action of $ G $ on the bundle which is compatible with the adjoint action on the base $ G $. \\
Let $ G $ be a locally compact quantum group and consider the $ G $-$ C^* $-algebra $ C^\red_0(G) $ with coaction $ \Delta $.
If $ G $ is regular the map $ \lambda: C^\red_0(G) \rightarrow
M(C^*_\red(G) \otimes C_0^\red(G)) $ given by
$$
\lambda(f) = \hat{W}^* (1 \otimes f) \hat{W}
$$
defines a continuous coaction. Moreover
\begin{align*}
(\ad(W) &\otimes \id)(\id \otimes \lambda) \Delta(f) = W_{12} \hat{W}_{23}^* W^*_{13}(1 \otimes 1 \otimes f) W_{13} \hat{W}_{23} W_{12}^* \\
&= \Sigma_{23} W_{13} W_{23} W^*_{12} \Sigma_{23}(1 \otimes 1 \otimes f) \Sigma_{23} W_{12} W_{23}^* W_{13}^* \Sigma_{23} \\
&= \Sigma_{23} W^*_{12} W_{23} \Sigma_{23}(1 \otimes 1 \otimes f) \Sigma_{23} W^*_{23} W_{12} \Sigma_{23}\\
&= W^*_{13} \hat{W}_{23}^*(1 \otimes 1 \otimes f) \hat{W}_{23} W_{13} = (\sigma \otimes \id) (\id \otimes \Delta) \lambda(f)
\end{align*}
shows that $ C^\red_0(G) $ together with $ \Delta $ and $ \lambda $ is a $ G $- $ \yd $-algebra.
More generally, we can consider a crossed product $ C^\red_0(G) \ltimes_\red A $ for a regular quantum group $ G $.
The dual coaction together with conjugation by $ \hat{W}^* $ as above yield a $ G $-$ \yd $-algebra structure on $ C^\red_0(G) \ltimes_\red A $. \\
There is another way to obtain a Yetter-Drinfeld-$ C^* $-algebra structure on a crossed product.
Let again $ G $ be a regular quantum group and let $ A $ be a $ G $-$ \yd $-algebra. We obtain a continuous coaction
$ \hat{\lambda}: C^*_\red(G)^\cop \ltimes_\red A \rightarrow M(C^*_\red(G) \otimes (C^*_\red(G)^\cop \ltimes_\red A)) $
by
$$
\hat{\lambda}(x) = \hat{W}_{12}^*(\id \otimes \lambda)(x)_{213} \hat{W}_{12}
$$
for $ x \in C^*_\red(G)^\cop \ltimes_\red A \subset \LH(\HH_G \otimes A) $. On the copy of $ A $ in the multiplier algebra of the crossed product
this coaction implements $ \lambda $, and on the copy of $ C^*_\red(G)^\cop = C^*_\red(G) $ it is given by the comultiplication $ \hat{\Delta} $
of $ C^*_\red(G) $.
In addition we have a continuous coaction $ \hat{\alpha}: C^*_\red(G)^\cop \ltimes_\red A \rightarrow M(C^\red_0(G) \otimes (C^*_\red(G)^\cop \ltimes_\red A))  $
given by
$$
\hat{\alpha}(x) = W_{12}^*(1 \otimes x)  W_{12}.
$$
Remark that on the copy of $ A $ in the multiplier algebra this coaction implements $ \alpha $, and
on the copy of $ C^*_\red(G)^\cop = C^*_\red(G) $ it implements the adjoint coaction. It is
straightforward to check that the crossed product $ C^*_\red(G)^\cop \ltimes_\red A $ becomes again a $ G $-$ \yd $-algebra
in this way. \\
The notion of a Yetter-Drinfeld-$ C^* $-algebra is closely related to coactions of the Drinfeld double.
Let us briefly recall the definition of the Drinfeld double in the context of locally compact quantum groups. It is
described as a special case of the double crossed product construction in \cite{BV}.
If $ G $ is a locally compact quantum group, then the reduced $ C^* $-algebra of functions on
the Drinfeld double $ \DD(G) $ is $ C^\red_0(\DD(G)) = C^\red_0(G) \otimes C^*_\red(G) $ with the comultiplication
$$
\Delta_{\DD(G)} = (\id \otimes \sigma \otimes \id)(\id \otimes \ad(W) \otimes \id)(\Delta \otimes \hat{\Delta}).
$$
This yields a locally compact quantum group $ \DD(G) $
which contains both $ G $ and $ \hat{G} $ as closed quantum subgroups. If $ G $ is regular then $ \DD(G) $ is again regular.
\begin{prop}
Let $ G $ be a locally compact quantum group and let $ \DD(G) $ be its Drinfeld double. Then a $ G $-Yetter-Drinfeld $ C^* $-algebra is
the same thing as a $ \DD(G) $-$ C^* $-algebra.
\end{prop}
\proof Let us first assume that $ A $ is a $ \DD(G) $-$ C^* $-algebra with coaction $ \gamma: A \rightarrow M(C^\red_0(\DD(G)) \otimes A) $.
Since $ G $ and $ \hat{G} $ are quantum subgroups of $ \DD(G) $ we obtain associated continuous
coactions $ \alpha: A \rightarrow M(C^\red_0(G) \otimes A) $ and
$ \lambda: A \rightarrow M(C^*_\red(G) \otimes A) $ by restriction. These coactions are determined by the conditions
$$
(\delta \otimes \id) \gamma = (\id \otimes \alpha)\gamma, \qquad (\hat{\delta} \otimes \id)\gamma = (\id \otimes \lambda) \gamma
$$
where the maps $ \delta: C^\red_0(\DD(G)) \rightarrow M(C^\red_0(\DD(G)) \otimes C^\red_0(G)) $ and
$ \hat{\delta}: C^\red_0(\DD(G)) \rightarrow M(C^\red_0(\DD(G)) \otimes C^*_\red(G)) $ are given by
$$
\delta = (\id \otimes \sigma) \ad(W_{23}) (\Delta \otimes \id), \qquad \hat{\delta} = \id \otimes \hat{\Delta}.
$$
We have
\begin{align*}
\ad(W_{23})(&\id \otimes \id \otimes \lambda)(\id \otimes \alpha)\gamma = \ad(W_{23}) (\id \otimes \id \otimes \lambda)(\delta \otimes \id)\gamma \\
&= \ad(W_{23}) (\delta \otimes \id \otimes \id)(\hat{\delta} \otimes \id)\gamma \\
&= \ad(W_{34}) (\id \otimes \sigma \otimes \id \otimes \id)\ad(W_{23}) (\Delta \otimes \hat{\Delta} \otimes \id)\gamma \\
&= (\id \otimes \sigma \otimes \id \otimes \id)\ad(W_{24} W_{23}) (\Delta \otimes \hat{\Delta} \otimes \id)\gamma \\
&= (\id \otimes \sigma \otimes \id \otimes \id)\ad((\id \otimes \hat{\Delta})(W)_{234}) (\Delta \otimes \hat{\Delta} \otimes \id)\gamma \\
&= (\id_{\DD(G)} \otimes \sigma \otimes \id)
(\id \otimes \hat{\Delta} \otimes \id \otimes \id)(\id \otimes \sigma \otimes \id) \ad(W_{23})(\Delta \otimes \id \otimes \id) \gamma \\
&= (\id \otimes \sigma \otimes \id)(\hat{\delta} \otimes \id \otimes \id)(\delta \otimes \id) \gamma \\
&= (\id \otimes \sigma \otimes \id)(\id \otimes \id \otimes \alpha)(\id \otimes \lambda)\gamma,
\end{align*}
and since the coaction $ \gamma $ is continuous this implies
$$
\ad(W_{12})(\id \otimes \lambda)\alpha = (\sigma \otimes \id)(\id \otimes \alpha) \lambda.
$$
It follows that we have obtained a $ G $-$ \yd $-algebra structure on $ A $. \\
Conversely, assume that $ A $ is equipped with a $ G $-$ \yd $-algebra structure. We define
a nondegenerate $ * $-homomorphism $ \gamma: A \rightarrow M(C^\red_0(\DD(G)) \otimes A) $ by
$$
\gamma = (\id \otimes \lambda)\alpha
$$
and compute
\begin{align*}
(\id \otimes \gamma)\gamma &= (\id \otimes \id \otimes \id \otimes \lambda)(\id \otimes \id \otimes \alpha)(\id \otimes \lambda)\alpha \\
&= (\id \otimes \id \otimes \id \otimes \lambda)(\id \otimes \sigma \otimes \id) \ad(W_{23}) (\id \otimes \id \otimes \lambda)(\id \otimes \alpha)\alpha \\
&= (\id \otimes \sigma \otimes \id \otimes \id)
\ad(W_{23})(\id \otimes \id \otimes \id \otimes \lambda)(\id \otimes \id \otimes \lambda)(\Delta \otimes \id)\alpha \\
&= (\id \otimes \sigma \otimes \id \otimes \id) \ad(W_{23}) (\Delta \otimes \hat{\Delta} \otimes \id)(\id \otimes \lambda)\alpha \\
&= (\id \otimes \sigma \otimes \id \otimes \id) \ad(W_{23}) (\Delta \otimes \hat{\Delta} \otimes \id)\gamma \\
&= (\Delta_{\DD(G)} \otimes \id)\gamma.
\end{align*}
It follows that $ \gamma $ is a continuous coaction which turns $ A $ into a
$ \DD(G) $-$ C^* $-algebra. \\
One checks easily that the two operations above are inverse to each other. \qed \\
We shall now define the braided tensor product $ A \twisted B $ of a $ G $-$ \yd $-algebra $ A $ with a $ G $-$ C^* $-algebra $ B $.
Observe first that the $ C^*$-algebra $ B $ acts on the Hilbert module $ \HH \otimes B $ by $ (\pi \otimes \id)\beta $ where
$ \pi: C^\red_0(G) \rightarrow \LH(\HH) $ denotes the defining representation on $ \HH = \HH_G $.
Similarly, the $ C^* $-algebra $ A $ acts on $ \HH \otimes A $ by
$ (\hat{\pi} \otimes \id)\lambda $ where $ \hat{\pi}: C^*_\red(G) \rightarrow \LH(\HH) $ is the
defining representation.
From this we obtain two $ * $-homomorphisms $ \iota_A = \lambda_{12}: A \rightarrow \LH(\HH \otimes A \otimes B) $ and
$ \iota_B = \beta_{13}: B \rightarrow \LH(\HH \otimes A \otimes B) $
by acting with the identity on the factor $ B $ and $ A $, respectively.
\begin{definition}
Let $ G $ be a locally compact quantum group, let $ A $ be a $ G $-$ \yd $-algebra and $ B $ a $ G $-$ C^* $-algebra.
With the notation as above, the braided tensor product $ A \twisted_G B $ is the
$ C^* $-subalgebra of $ \LH(\HH \otimes A \otimes B) $ generated by all elements
$ \iota_A(a) \iota_B(b) $ for $ a \in A $ and $ b \in B $.
\end{definition}
We will also write $ A \twisted B $ instead of $ A \twisted_G B $ if the quantum group $ G $ is clear from the context.
The braided tensor product $ A \twisted B $ is in fact equal to the closed linear span $ [\iota_A(A) \iota_B(B)] $.
This follows from proposition 8.3 in \cite{Vaesimprimitivity}, we
reproduce the argument for the convenience of the reader. Clearly it suffices to prove $ [\iota_A(A) \iota_B(B)] = [\iota_B(B) \iota_A(A)] $.
Using continuity of the coaction $ \lambda $ and $ \hat{V} = (J\hat{J} \otimes 1) W^* (\hat{J} J\otimes 1) $ we get
\begin{align*}
\lambda(A) &= [(\LH(\HH_G)_* \otimes \id \otimes \id) (\hat{\Delta} \otimes \id) \lambda(A)] \\
&= [(\LH(\HH_G)_* \otimes \id \otimes \id) (\hat{V}_{12} \lambda(A)_{13} \hat{V}_{12}^*)] \\
&= [(\LH(\HH_G)_* \otimes \id \otimes \id) (W^*_{12} \mu(A)_{13} W_{12})]
\end{align*}
where $ \mu(x) = (\hat{J} J \otimes 1) \lambda(x) (J \hat{J} \otimes 1) $ for $ x \in A $.
Since $ \beta: B \rightarrow M(C^\red_0(G) \otimes B) $ is a continuous coaction we obtain
\begin{align*}
[\lambda(A)_{12} \beta(B)_{13}] &= [(\LH(\HH_G)_* \otimes \id \otimes \id \otimes \id)(W^*_{12} \mu(A)_{13} W_{12} \beta(B)_{24})] \\
&= [(\LH(\HH_G)_* \otimes \id \otimes \id \otimes \id)(W^*_{12} \mu(A)_{13} W_{12} (\Delta \otimes \id)\beta(B)_{124})] \\
&= [(\LH(\HH_G)_* \otimes \id \otimes \id \otimes \id)(W^*_{12} \mu(A)_{13} \beta(B)_{24} W_{12})] \\
&= [(\LH(\HH_G)_* \otimes \id \otimes \id \otimes \id)((\Delta \otimes \id)\beta(B)_{124}) W^*_{12} \mu(A)_{13} W_{12})] \\
&= [(\LH(\HH_G)_* \otimes \id \otimes \id \otimes \id)(\beta(B)_{24}) W^*_{12} \mu(A)_{13} W_{12})] \\
&= [\beta(B)_{13}\lambda(A)_{12}]
\end{align*}
which yields the claim. It follows in particular that we have natural nondegenerate $ * $-homomorphisms
$ \iota_A: A \rightarrow M(A \twisted B) $ and $ \iota_B: B \rightarrow M(A \twisted B) $. \\
The braided tensor product $ A \twisted B $ becomes a $ G $-$ C^* $-algebra in a canonical way.
In fact, we have a nondegenerate
$ * $-homomorphism $ \alpha \boxtimes \beta: A \twisted B \rightarrow M(C^\red_0(G) \otimes (A \twisted B)) $
given by
\begin{align*}
(\alpha \twisted \beta)(\lambda(a)_{12} \beta(b)_{13}) &= W^*_{12} (\sigma \otimes \id)((\id \otimes \alpha)\lambda(a))_{123} \beta(b)_{24} W_{12} \\
&= (\id \otimes \lambda)\alpha(a)_{123} (\id \otimes \beta)\beta(b)_{124},
\end{align*}
and it is straightforward to check that $ \alpha \twisted \beta $ defines a continuous coaction of $ C^{\red}_0(G) $ such that the $ * $-homomorphisms
$ \iota_A $ and $ \iota_B $ are $ G $-equivariant. \\
If $ B $ is a $ G $-$ \yd $-algebra with coaction $ \gamma: B \rightarrow M(C^*_\red(G) \otimes B) $
then we obtain a nondegenerate
$ * $-homomorphism $ \lambda \twisted \gamma: A \twisted B \rightarrow M(C^*_\red(G) \otimes (A \twisted B)) $
by the formula
\begin{align*}
(\lambda \twisted \gamma)(\lambda(a)_{12} \beta(b)_{13}) &= \hat{W}_{12}^* \lambda(a)_{23}(\sigma \otimes \id)((\id \otimes \gamma)\beta(b))_{124} \hat{W}_{12} \\
&= (\id \otimes \lambda)\lambda(a)_{123} (\id \otimes \beta)\gamma(b)_{124}.
\end{align*}
In the same way as above one finds that $ \lambda \twisted \gamma $ yields a continuous coaction of $ C^*_\red(G) $
such that $ \iota_A $ and $ \iota_B $ are $ \hat{G} $-equivariant. From
the equivariance of $ \iota_A $ and $ \iota_B $ it follows that $ A \twisted B $
together with the coactions $ \alpha \twisted \beta $ and $ \lambda \twisted \gamma $
becomes a $ G $-$ \yd $-algebra. \\
If $ A $ is a $ G $-$ \yd $-algebra and $ f: B \rightarrow C $ a possibly degenerate equivariant $ * $-homomorphism of $ G $-$ C^* $-algebras,
then we obtain an induced $ * $-homomorphism $ M_\KH(\KH \otimes A \otimes B)
\rightarrow M_\KH(\KH \otimes A \otimes C) $ between the relative multiplier algebras. Since
$ A \twisted B \subset M(\KH \otimes A \otimes B) $
is in fact contained in $ M_\KH(\KH \otimes A \otimes B) $, this map restricts to an equivariant $ * $-homomorphism
$ \id \twisted f: A \twisted B \rightarrow A \twisted C $. It follows that the braided tensor product defines a functor
$ A \twisted - $ from $ G \Alg $ to $ G \Alg $.
Similarly, if $ f: A \rightarrow B $ is a homomorphism of $ G $-$ \yd $-algebras we obtain for
every $ G $-algebra $ C $ an equivariant $ * $-homomorphism $ f \twisted \id: A \twisted C \rightarrow B \twisted C $
and a functor $ - \twisted C $ from $ \DD(G)\Alg $ to $ G \Alg $.
There are analogous functors $ A \twisted - $ and $ - \twisted C $ from $ \DD(G)\Alg $ to $ \DD(G)\Alg $
if we consider $ G $-$ \yd $-algebras in the second variable. \\
Assume now that $ A $ and $ B $ are $ G $-$ \yd $-algebras and that $ C $ is a $ G $-$ C^* $-algebra. According to
our previous observations we can form the braided tensor products $ (A \twisted B) \twisted C $ and $ A \twisted (B \twisted C) $,
respectively. We have
\begin{align*}
(A \twisted B) \twisted C &= [(\id \otimes \lambda^A) \lambda^A(A)_{123} (\id \otimes \beta)\lambda^B(B)_{124} \gamma(C)_{15}] \\
&= [\hat{W}_{12}^* \lambda^A(A)_{23} \hat{W}_{12} (\id \otimes \beta)\lambda^B(B)_{124} \hat{W}^*_{12} \Sigma_{12} W^*_{12} \Sigma_{12}
\gamma(C)_{15}] \\
&= [\hat{W}_{12}^* \lambda^A(A)_{23} \Sigma_{12} (\id \otimes \lambda^B)\beta(B)_{124} W^*_{12}
\Sigma_{12} \gamma(C)_{15}\Sigma_{12} W_{12} \Sigma_{12} \hat{W}_{12}] \\
&= [\hat{W}_{12}^* \Sigma_{12} \lambda^A(A)_{13} (\id \otimes \lambda^B)\beta(B)_{124} (\id \otimes \gamma)\gamma(C)_{125} \Sigma_{12} \hat{W}_{12}] \\
&\cong [\lambda^A(A)_{13} (\beta \twisted \gamma)(B \twisted C)_{1245}] \cong A \twisted (B \twisted C),
\end{align*}
and the resulting isomorphism $ (A \twisted B) \twisted C \cong A \twisted (B \twisted C) $ is $ G $-equivariant.
If $ C $ is a $ G $-$ \yd $-algebra then this isomorphism is in addition $ \hat{G} $-equivariant.
We conclude that the braided tensor product is associative in a natural way. \\
If $ B $ is a trivial $ G $-algebra then the braided tensor product $ A \twisted B $ is isomorphic
to $ A \otimes B $ with the coaction induced from $ A $. Similarly, if the coaction of $ C^*_\red(G) $ on the $ G $-$ \yd $-algebra
$ A $ is trivial then $ A \twisted B $ is isomorphic to $ A \otimes B $.
Recall that if $ G $ is a locally compact group we may view all $ G $-algebras as $ G $-$ \yd $-algebras
with the trivial coaction of the group $ C^* $-algebra. In this case the braided tensor product reduces to the ordinary tensor product
of $ G $-$ C^* $-algebras with the diagonal $ G $-action. For general quantum groups the braided tensor product should be viewed as
a substitute for the latter construction. \\
Following an idea of Vaes, we shall now discuss the compatibility of the braided tensor product with induction and restriction.
Let $ G $ be a strongly regular quantum group and let $ H \subset G $ be a closed quantum subgroup determined by
the faithful normal $ * $-homomorphism $ \hat{\pi}: \mathcal{L}(H) \rightarrow \mathcal{L}(G) $.
Keeping our notation from section \ref{secqg}, we denote by $ I $ the corresponding von-Neumann algebraic imprimitivity bimodule
for $ \mathcal{L}(G)^\cop \ltimes L^\infty(G/H) $ and $ \mathcal{L}(H)^\cop $, and by $ \I \subset I $ the $ C^* $-algebraic imprimitivity bimodule
for $ C^*_\red(G)^\cop \ltimes_\red C^\red_0(G/H) $ and $ C^*_\red(H)^\cop $.
\begin{prop} \label{indyd}
Let $ G $ be a strongly regular quantum group and let $ H \subset G $ be a closed quantum subgroup.
If $ A $ is an $ H $-$ \yd $-algebra then the induced $ C^* $-algebra $ \ind_H^G(A) $ is a $ G $-$ \yd $-algebra
in a natural way.
\end{prop}
\proof Let $ \alpha: A \rightarrow M(C^\red_0(H) \otimes A) $ be the coaction of $ C^\red_0(H) $ on $ A $.
From the construction of $ \ind_H^G(A) $ in \cite{Vaesimprimitivity} we have the induced coaction
$ \ind(\alpha): \ind_H^G(A) \rightarrow M(C^\red_0(G) \otimes A) $
given by $ \ind(\alpha)(x) = (W^*_G)_{12}(1 \otimes x) (W_G)_{12} $ for
$ x \in \ind_H^G(A) \subset \LH(\HH_G \otimes A) $. Our task is to define a continuous coaction of $ C^*_\red(G) $ on $ \ind_H^G(A) $
satisfying the $ \yd $-condition. \\
Denote by $ \lambda: A \rightarrow M(C^*_\red(H) \otimes A) $
the coaction which determines the $ H $-$ \yd $-algebra structure on $ A $. This coaction
induces a coaction $ \res(\lambda): A \rightarrow M(C^*_\red(G) \otimes A) $ because
$ H \subset G $ is a closed quantum subgroup.
Since $ A $ is an $ H $-$ \yd $-algebra we have in addition the coaction $ \hat{\lambda} $ of
$ C^*_\red(H) $ on the crossed product $ C^*_\red(H)^\cop \ltimes_\red A $,
and a corresponding coaction $ \res(\hat{\lambda}) $ of $ C^*_\red(G) $. \\
We abbreviate $ B = C^*_\red(H)^\cop \ltimes_\red A $ and consider the Hilbert $ B $-module $ \E = B $ with
the corepresentation $ X = W_H \otimes \id \in M(C^\red_0(H) \otimes \KH(\E)) = M(C^\red_0(H) \otimes B)$.
The corresponding induced Hilbert $ B $-module $ \ind_H^G(\E) $ is constructed in \cite{Vaesimprimitivity} such that
$$
\HH_G \otimes \ind_H^G(\E) \cong I \otimes_{\pi_l} \F
$$
where $ \F = \HH_G \otimes \E $ and the strict $ * $-homomorphism $ \pi_l: \mathcal{L}(H) \rightarrow \LH(\F) $ is determined by
$ (\id \otimes \pi_l)(W_H) = (\id \otimes \hat{\pi})(W_H)_{12} X_{13} $. \\
Let us define a coaction on $ I \otimes_{\pi_l} \F $ as follows. On $ I $ we have the adjoint action $ \eta: I \rightarrow \mathcal{L}(G) \otimes I $ given by
$$
\eta(v) = \hat{W}^*_G(1 \otimes v) (\hat{\pi} \otimes \id)(\hat{W}_H)
$$
which is compatible with the coaction $ (\hat{\pi} \otimes \id)\hat{\Delta}_H: \mathcal{L}(H) \rightarrow \mathcal{L}(G) \otimes \mathcal{L}(H) $.
In addition consider the coaction $ \beta_\F $ of $ C^*_\red(G) $ on $ \F $
given by $ \beta_\F = (\sigma \otimes \id)(\id \otimes \res(\hat{\lambda})) $.
By construction, $ \beta_\F $ is compatible with the coaction $ \res(\hat{\lambda}) $ on $ B $.
Moreover, the $ * $-homomorphism $ \pi_l: \mathcal{L}(H) \rightarrow \LH(\F) $
is covariant in the sense that
$$
(\id \otimes \pi_l)(\hat{\pi} \otimes \id)\hat{\Delta}_H(x) =
\ad_{\beta_\F}(\pi_l(x))
$$
in $ \LH(\HH_G \otimes \F) $ for all $ x \in \mathcal{L}(H) $. According to proposition 12.13 in \cite{Vaesimprimitivity} we obtain a product coaction
of $ C^*_\red(G) $ on $ I \otimes_{\pi_l} \F $. \\
Under the above isomorphism, this product coaction leaves invariant the natural representations of $ L^\infty(G)' $ and $ \mathcal{L}(G)' $ on the
first tensor factor of the left hand side.
Hence there is an induced coaction
$ \gamma: \ind_H^G(\E) \rightarrow M(C^*_\red(G) \otimes \ind_H^G(\E)) $
on $ \ind_H^G(\E) $. Using the identification
$ \ind_H^G(\E) \cong [(\I \otimes 1) \alpha(A)] $ we see that $ \gamma $ is given by
$$
\gamma((v \otimes 1)\alpha(a)) = (\eta(v) \otimes 1) (\id \otimes \alpha)\res(\lambda)(a)
$$
for $ v \in \I $ and $ a \in A $. Since $ \KH(\ind_H^G(\E)) = C^*_\red(G)^\cop \ltimes_\red \ind_H^G(A) $ we obtain
a coaction $ \ad_\gamma $ on $ C^*_\red(G)^\cop \ltimes_\red \ind_H^G(A) $.
By construction, the coaction $ \ad_\gamma $ commutes with the dual coaction and is given
by $ \hat{\Delta} $ on the copy of $ C^*_\red(G)^\cop $.
It follows that $ \ad_\gamma $ induces a
continuous coaction $ \delta: \ind_H^G(A) \rightarrow M(C^*_\red(G) \otimes \ind_H^G(A)) $. Explicitly, this coaction is given by
$$
\delta(x) = (\hat{W}_G^*)_{12} (\sigma \otimes \id)(\id \otimes \res(\lambda))(x) (\hat{W}_G)_{12}
$$
for $ x \in \ind_H^G(A) \subset \LH(\HH_G \otimes A) $. Writing $ W_G = W $ and $ \hat{W}_G = \hat{W} $ we calculate
\begin{align*}
\ad&(W_{12})(\id \otimes \delta)\ind(\alpha)(x) \\
&= W_{12} \hat{W}^*_{23} \Sigma_{23}(\id \otimes \id \otimes \res(\lambda)) (W^*_{12}(1 \otimes x) W_{12}) \Sigma_{23} \hat{W}_{23} W_{12}^* \\
&= \Sigma_{23} W_{13} W_{23} W^*_{12} (\id \otimes \id \otimes \res(\lambda))(1 \otimes x) W_{12} W_{23}^* W_{13}^* \Sigma_{23} \\
&= \Sigma_{23} W^*_{12} W_{23} (\id \otimes \id \otimes \res(\lambda))(1 \otimes x) W^*_{23} W_{12} \Sigma_{23}\\
&= W^*_{13} \hat{W}^*_{23} \Sigma_{23}(\id \otimes \id \otimes \res(\lambda))(1 \otimes x) \Sigma_{23} \hat{W}_{23} W_{13} \\
&= (\sigma \otimes \id) (\id \otimes \ind(\alpha)) \delta(x)
\end{align*}
which shows that $ \ind(\alpha) $ and $ \delta $ combine to turn $ \ind_H^G(A) $ into a $ G $-$ \yd $-algebra. \qed \\
Let $ G $ be a regular locally compact quantum group and let $ A $ be a $ G $-$ \yd $-algebra with coactions $ \alpha $ and $ \lambda $.
As explained above, the crossed product $ C^*_\red(G)^\cop \ltimes_\red A $ is again a $ G $-$ \yd $-algebra in a natural way.
Moreover let $ B $ be a $ G $-algebra with coaction $ \beta $ and observe
\begin{align*}
C^*_\red(G)^\cop \ltimes_\red (A \twisted_G B) &=
[(C^*_\red(G) \otimes \id \otimes \id \otimes \id)(\id \otimes \lambda)\alpha(A)_{123} (\id \otimes \beta)\beta(B)_{124}] \\
&\cong [\hat{\Delta}(C^*_\red(G))_{21} W_{12} (\id \otimes \lambda)\alpha(A)_{123} W^*_{12} \beta(B)_{24}] \\
&= [\hat{\Delta}(C^*_\red(G))_{21} (\id \otimes \alpha)\lambda(A)_{213} \beta(B)_{24}] \\
&\cong [\hat{\lambda}(C^*_\red(G)^\cop \ltimes_\red A)_{12} \beta(B)_{13}] = (C^*_\red(G)^\cop \ltimes_\red A) \twisted_G B.
\end{align*}
Under this isomorphism the dual coaction on the left hand side corresponds to the coaction determined by the dual coaction on
$ C^*_\red(G)^\cop \ltimes_\red A $ and the trivial coaction on $ B $ on the right hand side.
As a consequence we obtain the following lemma.
\begin{lemma} \label{indtwistedlemma}
Let $ G $ be a regular locally compact quantum group, let $ A $ be a $ G $-$ \yd $-algebra and let $ B $ be a $ G $-algebra. Then
there is a natural $ C^*_\red(G)^\cop $-colinear isomorphism
$$
C^*_\red(G)^\cop \ltimes_\red (A \twisted_G B) \cong (C^*_\red(G)^\cop \ltimes_\red A) \twisted_G B.
$$
\end{lemma}
After these preparations we shall now describe the compatibility of restriction, induction and braided tensor products.
\begin{theorem} \label{citensor}
Let $ G $ be a strongly regular quantum group and let $ H \subset G $ be a closed quantum subgroup.
Moreover let $ A $ be an $ H $-$ \yd $-algebra and let $ B $ be a $ G $-algebra.
Then there is a natural $ G $-equivariant isomorphism
$$
\ind_H^G(A \twisted_H \res^G_H(B)) \cong \ind_H^G(A) \twisted_G B.
$$
\end{theorem}
\proof Note that the case $ A = \mathbb{C} $ with the trivial action is treated in \cite{Vaesimprimitivity}.
We denote by $ \res(\beta) $ the restriction to $ C^\red_0(H) $ of the coaction $ \beta: B \rightarrow M(C^\red_0(G) \otimes B) $.
Moreover let $ \res(\lambda) $ be the push-forward of the
coaction $ \lambda: A \rightarrow M(C^*_\red(H) \otimes A) $ to $ C^*_\red(G) $. Then
\begin{align*}
[(\id \otimes &\id \otimes \beta)(\lambda(A)_{12} \res(\beta)(B)_{13})] = [\lambda(A)_{12} (\id \otimes \id \otimes \beta)\res(\beta)(B)_{134}] \\
&= [\lambda(A)_{12} (\id \otimes \hat{\pi})(W_H^*)_{13} \beta(B)_{34} (\id \otimes \hat{\pi})(W_H)_{13}] \\
&= [(\id \otimes \hat{\pi})(W_H^*)_{13} ((\hat{\pi} \otimes \id \otimes \id) (\hat{\Delta}_H \otimes \id)\lambda(A))_{312}
\beta(B)_{34} (\id \otimes \hat{\pi})(W_H)_{13}] \\
&\cong [((\hat{\pi} \otimes \id \otimes \id) (\id \otimes \lambda)\lambda(A))_{312} \beta(B)_{34}] \\
&= [(\lambda \otimes \id \otimes \id)(\res(\lambda)(A)_{21} \beta(B)_{23})],
\end{align*}
and hence
$$
A \boxtimes_H \res^G_H(B) \cong [\res(\lambda)(A)_{12} \beta(B)_{13}].
$$
Writing $ W_G = W $ we conclude
\begin{align*}
&C^*_\red(G)^{\cop} \ltimes_\red \ind_H^G(A \twisted_H \res^G_H(B)) \\
&\cong [(\I \otimes \id \otimes \id \otimes \id)(\id \otimes \res(\lambda))\alpha(A)_{123}
(\id \otimes \beta)\res(\beta)(B)_{124}(\I^* \otimes \id \otimes \id \otimes \id)] \\
&= [((\id \otimes \res(\lambda))(\I \otimes \id)\alpha(A)(\I^* \otimes \id))_{123}
(\Delta_G \otimes \id)\beta(B)_{124}] \\
&= [W^*_{12} ((W \otimes \id)(\id \otimes \res(\lambda))((\I \otimes \id)\alpha(A)(\I^* \otimes \id)) (W^* \otimes 1))_{123} \beta(B)_{24} W_{12}]
\end{align*}
using that $ [(\I \otimes \id)\res(\beta)(B)] = [\beta(B)(\I \otimes \id)] $ for the restricted coaction $ \res(\beta) $.
Moreover
\begin{align*}
&[W^*_{12} ((W \otimes \id)(\id \otimes \res(\lambda))((\I \otimes \id)\alpha(A)(\I^* \otimes \id)) (W^* \otimes 1))_{123} \beta(B)_{24} W_{12}] \\
&\cong [((W \otimes \id)(\id \otimes \res(\lambda))((\I \otimes \id)\alpha(A)(\I^* \otimes \id)) (W^* \otimes 1))_{123} \beta(B)_{24}] \\
&= [\hat{\delta}(C^*_\red(G)^\cop \ltimes_\red \ind_H^G(A))_{213} \beta(B)_{24}] \\
&\cong (C^*_\red(G)^\cop \ltimes_\red \ind_H^G(A)) \twisted_G B
\end{align*}
where $ \hat{\delta}: C^*_\red(G)^{\cop} \ltimes_\red \ind_H^G(A) \rightarrow M(C^*_\red(G) \otimes (C^*_\red(G)^{\cop} \ltimes_\red \ind_H^G(A)) $
is the natural coaction on the crossed product of the $ G $-$ \yd $-algebra $ \ind_H^G(A) $. \\
Under these identifications, the dual coaction on $ C^*_\red(G)^{\cop} \ltimes_\red \ind_H^G(A \twisted_H \res^G_H(B)) $ corresponds on
$ (C^*_\red(G)^\cop \ltimes_\red \ind_H^G(A)) \twisted_G B $ to the dual
coaction on the crossed product and the trivial coaction on $ B $. As a consequence,
using lemma \ref{indtwistedlemma} we obtain a $ C^*_\red(G)^\cop $-colinear isomorphism
\begin{align*}
C^*_\red(G)^\cop \ltimes_\red \ind_H^G(A \twisted_H \res^G_H(B)) \cong C^*_\red(G)^\cop \ltimes_\red (\ind_H^G(A) \twisted_G B).
\end{align*}
Moreover, the element $ W \otimes \id \in M(C^\red_0(G) \otimes C^*_\red(G)^\cop \ltimes_\red \ind_H^G(A \twisted_H \res^G_H(B))) $
is mapped to $ W \otimes \id \in M(C^\red_0(G) \otimes C^*_\red(G)^\cop \ltimes_\red (\ind_H^G(A) \twisted_G B)) $ under this isomorphism.
Due to theorem 6.7 in \cite{Vaesimprimitivity} this shows that there is a $ G $-equivariant isomorphism
$$
\ind_H^G(A \twisted_H \res^G_H(B)) \cong \ind_H^G(A) \twisted_G B
$$
as desired. \qed \\
We also need braided tensor products of Hilbert modules. Since the constructions and arguments are similar to the algebra case treated above our
discussion will be rather brief. Assume that $ A $ is a $ G $-$ \yd $-algebra and that $ B $ is a $ G $-algebra.
Moreover let $ \E_A $ be a $ \DD(G) $-Hilbert module and let $ \F_B $ be a $ G $-Hilbert module.
As in the algebra case, a $ \DD(G) $-Hilbert module $ \E $ is the same thing as a Hilbert module equipped with continuous coactions $ \alpha_\E $ of
$ C^\red_0(G) $ and $ \lambda_\E $ of $ C^*_\red(G) $ satisfying the Yetter-Drinfeld compatibility condition in the sense that
$$
(\sigma \otimes \id)(\id \otimes \alpha_\E)\lambda_\E = (\ad(W) \otimes \id) (\id \otimes \lambda_\E)\alpha_\E
$$
where $ \ad(W) $ is the adjoint action. \\
The braided tensor product of $ \E $ and $ \F $ is defined as
$$
\E \twisted_G \F = [\lambda_\E(\E)_{12} \beta_\F(\F)_{13}] \subset M_\KH(\KH \otimes \E \otimes \F)
$$
where $ \lambda_\E $ denotes the coaction of $ C^*_\red(G) $ on $ \E $ and $ \beta_\F $ is the coaction of $ C^\red_0(G) $ on $ \F $.
One has $ [\lambda_\E(\E)_{12} \beta_\F(\F)_{13}] = [\beta_\F(\F)_{13} \lambda_\E(\E)_{12}] $,
and $ \E \twisted_G \F $ is closed under right multiplication by elements from $ A \twisted_G B \subset M_\KH(\KH \otimes A \otimes B) $.
Moreover the restriction to $ \E \twisted_G \F $ of the scalar product of $ M_\KH(\KH \otimes \E \otimes \F) $ takes values in $ A \twisted_G B $.
It follows that $ \E \twisted_G \F $ is a Hilbert-$ A \twisted_G B $-module. \\
As in the algebra case there is a continuous coaction of $ C^\red_0(G) $ on $ \E \twisted_G \F $
given by
$$
\ad(W^*_{12}) (\sigma \otimes \id)(\id \otimes \alpha_\E \otimes \id).
$$
Similarly, if $ B $ is a $ G $-$ \yd $-algebra and $ \F $ is a $ \DD(G) $-Hilbert module we have a continuous
$ C^*_\red(G) $-coaction. The braided tensor product becomes a $ \DD(G) $-Hilbert module in this case. \\
There are canonical nondegenerate $ * $-homomorphisms $ \KH(\E) \rightarrow \LH(\E \twisted_G \F) $ and $ \KH(\F) \rightarrow \LH(\E \twisted_G \F) $.
Combining these homomorphisms yields an identification $ \KH(\E) \twisted_G \KH(\F) \cong \KH(\E \twisted_G \F) $. \\
We conclude this section with a discussion of stability properties.
\begin{prop}\label{twistedstability}
Let $ G $ be a regular locally compact quantum group and let $ A $ be a $ G $-$ \yd $-algebra.
\begin{bnum}
\item[a)] For every $ G $-$ C^* $-algebra $ B $ there is a natural $ G $-equivariant
Morita equivalence
$$
(\KH_{\DD(G)} \otimes A) \twisted_G B \sim_M A \twisted_G B.
$$
If $ B $ is a $ G $-$ \yd $-algebra this Morita equivalence is $ \DD(G) $-equivariant.
\item[b)] For every $ G $-$ C^* $-algebra $ B $ there is a natural $ G $-equivariant Morita equivalence
$$
A \twisted_G (\KH_G \otimes B) \sim_M A \twisted_G B.
$$
If $ B $ is a $ G $-$ \yd $-algebra there is a natural $ \DD(G) $-equivariant Morita equivalence
$$
A \twisted_G (\KH_{\DD(G)} \otimes B) \sim_M A \twisted_G B.
$$
\end{bnum}
\end{prop}
\proof We consider the coaction of $ \DD(G) $ on $ \HH_{\DD(G)} $ coming from the regular representation.
From \cite{BV} we know that the corresponding corepresentation of $ C^\red_0(G) $ on $ \HH_{\DD(G)} = \HH_G \otimes \HH_G $
is $ W_{12} \in M(C^\red_0(G) \otimes \KH_{\DD(G)}) $. The corresponding corepresentation
of $ C^*_\red(G) $ on $ \HH_{\DD(G)} $ is given by $ Z^*_{23} \hat{W}_{13} Z_{23} $ where
$ Z = W(J \otimes \hat{J}) W(J \otimes \hat{J}) $. \\
To prove $ a) $ we observe that $ Z^*_{23} \hat{W}_{13} Z_{23} $ implements a $ G $-equivariant isomorphism
\begin{align*}
(\HH_{\DD(G)} \otimes A) \twisted_G B &= [Z^*_{23} \hat{W}^*_{13} Z_{23} \sigma_{12}(\id \otimes \lambda)(\HH_{\DD(G)} \otimes A)_{123} \beta(B)_{14}] \\
&\cong \HH_{\DD(G)} \otimes [\lambda(A)_{12} \beta(B)_{13}] = \HH_{\DD(G)} \otimes (A \twisted_G B)
\end{align*}
of Hilbert modules. This yields
\begin{align*}
(\KH_{\DD(G)} &\otimes A) \twisted_G B \cong \KH((\HH_{\DD(G)} \otimes A) \twisted_G B) \\
&\cong \KH(\HH_{\DD(G)} \otimes (A \twisted_G B)) \sim_M \KH(A \twisted_G B) = A \twisted_G B
\end{align*}
in a way compatible with the coaction of $ C^\red_0(G) $.
If $ B $ is a $ G $-$ \yd $-algebra the above isomorphisms and the Morita equivalence are $ \DD(G) $-equivariant. The assertions in $ b) $
are proved in a similar fashion. \qed

\section{The equivariant Kasparov category} \label{seckkg}

In this section we first review the definition of equivariant Kasparov theory given by Baaj and Skandalis \cite{BSKK}.
Then we explain how to extend several standard results from the case of locally compact groups to the setting of regular locally compact quantum groups.
In particular, we adapt the Cuntz picture of $ KK $-theory \cite{Cuntznewlook}
to show that equivariant $ KK $-classes can be described by homotopy classes of
equivariant homomorphisms. As a consequence, we obtain the universal property of equivariant Kasparov theory. We describe its structure as a triangulated category
and discuss the restriction and induction functors. Finally, based on the construction of the braided tensor product in the
previous section we construct exterior products in equivariant $ KK $-theory. \\
Let us recall the definition of equivariant Kasparov theory \cite{BSKK}. For simplicity we will assume that all $ C^* $-algebras are separable.
Let $ S $ be a Hopf-$ C^* $-algebra and let $ A $ and $ B $ be graded $ S $-$ C^* $-algebras. An $ S $-equivariant Kasparov $ A $-$ B $-module
is a countably generated graded $ S $-equivariant Hilbert $ B $-module $ \E $ together with an $ S $-colinear graded $ * $-homomorphism
$ \phi: A \rightarrow \LH(\E) $ and an odd operator $ F \in \LH(\E) $ such that
$$
[F, \phi(a)], \qquad (F^2 - 1) \phi(a), \qquad (F - F^*)\phi(a)
$$
are contained in $ \KH(\E) $ for all $ a \in A $ and $ F $ is almost invariant in the sense that
$$
(\id \otimes \phi)(x)(1 \otimes F - \ad_\lambda(F)) \subset S \otimes \KH(\E)
$$
for all $ x \in S \otimes A $. Here $ S \otimes \KH(\E) = \KH(S \otimes \E) $ is viewed as a subset of $ \LH(S \otimes \E) $
and $ \ad_\lambda $ is the adjoint coaction associated to the given coaction $ \lambda: \E \rightarrow M(S \otimes \E) $ on $ \E $.
Two $ S $-equivariant Kasparov $ A $-$ B $-modules $ (\E_0, \phi_0, F_0) $ and $ (\E_1, \phi_1, F_1) $ are called unitarly equivalent if there is
an $ S $-colinear unitary $ U \in \LH(\E_0, \E_1) $ of degree zero such that $ U \phi_0(a) = \phi_1(a) U $ for all $ a \in A $
and $ F_1 U = U F_0 $. We write $ (\E_0, \phi_0, F_0) \cong (\E_1, \phi_1, F_1) $ in this case.
Let $ E_S(A,B) $ be the set of unitary equivalence classes of $ S $-equivariant Kasparov $ A $-$ B $-modules.
This set is functorial for graded $ S $-colinear $ * $-homomorphisms
in both variables. If $ f: B_1 \rightarrow B_2 $ is a graded $ S $-colinear $ * $-homomorphism and $ (\E, \phi, F) $ is an
$ S $-equivariant Kasparov $ A $-$ B_1 $-module, then
$$
f_*(\E, \phi, F) = (\E \cotimes_f B_2, \phi \cotimes \id, F \cotimes 1)
$$
is the corresponding Kasparov $ A $-$ B_2 $-module.
A homotopy between $ S $-equivariant Kasparov $ A $-$ B $-modules $ (\E_0, \phi_0, F_0) $ and $ (\E_1, \phi_1, F_1) $ is an
$ S $-equivariant Kasparov $ A $-$ B[0,1] $-module
$ (\E, \phi, F) $ such that $ (\ev_t)_*(\E, \phi, F) \cong (\E_t, \phi_t, F_t) $ for $ t = 0,1 $. Here
$ B[0,1] = B \otimes C[0,1] $ where $ C[0,1] $ is equipped with the trivial action and grading
and $ \ev_t: B[0,1] \rightarrow B $ is evaluation at $ t $.
\begin{definition}
Let $ S $ be a Hopf-$ C^* $-algebra and let $ A $ and $ B $ be graded $ S $-$ C^* $-algebras. The $ S $-equivariant Kasparov group $ KK^S(A,B) $ is the set
of homotopy classes of $ S $-equivariant Kasparov $ A $-$ B $-modules.
\end{definition}
In the definition of $ KK^S(A,B) $ one can restrict to Kasparov triples $ (\E, \phi, F) $ which are
essential in the sense that $ [\phi(A) \E] = \E $, compare \cite{Meyerkkg}.
We note that $ KK^S(A,B) $ becomes an abelian group with addition given by the direct sum of Kasparov modules.
Many properties of ordinary $ KK $-theory carry over to the $ S $-equivariant situation, in particular the construction of the
Kasparov composition product and Bott periodicity \cite{BSKK}.
As usual we write $ KK^S_0(A,B) = KK^S(A,B) $ and let $ KK^S_1(A,B) $ be the odd $ KK $-group obtained by suspension in either variable.
In the case $ S = C_0(G) $ for a locally compact group $ G $ one reobtains the definition of $ G $-equivariant
$ KK $-theory \cite{Kasparov2}. \\
Our first aim is to establish the Cuntz picture of equivariant $ KK $-theory in the setting of regular locally compact quantum groups.
This can be done parallel to the account in the group case given by Meyer \cite{Meyerkkg}. For convenience we restrict ourselves to
trivially graded $ C^* $-algebras and present a short argument using Baaj-Skandalis duality. \\
Let $ S $ be a Hopf-$ C^* $-algebra and let $ A_1 $ and $ A_2 $ be $ S $-$ C^* $-algebras. Consider the free
product $ A_1 * A_2 $ together with the canonical $ * $-homomorphisms $ \iota_j: A_j \rightarrow A_1 * A_2 $ for $ j = 1,2 $.
We compose the coaction $ \alpha_j: A_j \rightarrow M_S(S \otimes A_j) $ with the $ * $-homomorphism
$ M_S(S \otimes A_j) \rightarrow M_S(S \otimes (A_1 * A_2)) $ induced by $ \iota_j $ and combine these maps to
obtain a $ * $-homomorphism $ \alpha: A_1 * A_2 \rightarrow M_S(S \otimes (A_1 * A_2)) $. This map
satisfies all properties of a continuous coaction in the sense of definition \ref{defcoaction} except that it is not
obvious wether $ \alpha $ is always injective. If necessary, this technicality can be overcome by passing to a quotient
of $ A_1 * A_2 $. More precisely, on $ A_1 *^S A_2 = (A_1 * A_2)/\ker(\alpha) $ the map $ \alpha $ induces
the structure of an $ S $-$ C^* $-algebra, and we have canonical $ S $-colinear $ * $-homomorphisms $ A_j \rightarrow A_1 *^S A_2 $
for $ j = 1,2 $ again denoted by $ \iota_j $.
The resulting $ S $-$ C^* $-algebra is universal for pairs of $ S $-colinear $ * $-homomorphisms $ f_1: A_1 \rightarrow C $ and
$ f_2 : A_2 \rightarrow C $ into $ S $-$ C^* $-algebras $ C $. That is, for any such pair of $ * $-homomorphisms there exists a unique
$ S $-colinear $ * $-homomorphism $ f: A_1 *^S A_2 \rightarrow C $ such that $ f \iota_j = f_j $ for $ j = 1,2 $.
By abuse of notation, we will still write $ A_1 * A_2 $ instead of $ A_1 *^S A_2 $ in the sequel.
We point out that in the arguments below we could equally well work with the ordinary free product together with its possibly noninjective coaction. \\
Let $ A $ be an $ S $-$ C^* $-algebra and consider $ QA = A * A $. The algebra $ \KH \otimes QA $ is
$ S $-colinearly homotopy equivalent to $ \KH \otimes (A \oplus A) $ where $ \KH $ denotes the algebra of compact operators on a separable
Hilbert space $ \HH $. Moreover there is an extension
\begin{equation*}
\xymatrix{
0 \ar@{->}[r] & qA \ar@{->}[r] & QA \ar@{->}[r]^{\pi} & A \ar@{->}[r] & 0
}
\end{equation*}
of $ S $-$ C^* $-algebras with $ S $-colinear splitting, here $ \pi $ is the homomorphism associated to the pair $ f_1 = f_2 = \id_A $ and
$ qA $ its kernel. \\
We shall now restrict attention from general Hopf-$ C^* $-algebras to regular locally compact quantum groups and state the
Baaj-Skandalis duality theorem \cite{BSKK}, \cite{BSUM}.
\begin{theorem}\label{BSduality}
Let $ G $ be a regular locally compact quantum group and let $ S = C^\red_0(G) $ and $ \hat{S} = C^*_\red(G)^\cop $.
For all $ S $-$ C^* $-algebras $ A $ and $ B $ there is a canonical isomorphism
$$
J_S: KK^S(A,B) \rightarrow KK^{\hat{S}}(\hat{S} \ltimes_\red A, \hat{S} \ltimes_\red B)
$$
which is multiplicative with respect to the composition product.
\end{theorem}
For our purposes it is important that under this isomorphism the class of an $ S $-equivariant Kasparov $ A $-$ B $-module $ (\E, \phi, F) $ is mapped
to the class of an $ \hat{S} $-equivariant Kasparov module $ (J_S(\E), J_S(\phi), J_S(F)) $ with an operator $ J_S(F) $ which is exactly invariant under the
coaction of $ \hat{S} $. \\
Let $ G $ be a regular locally compact quantum group and let $ \E $ and $ \F $ be $ G $-Hilbert $ B $-modules which are isomorphic as Hilbert $ B $-modules.
Then we have a $ G $-equivariant isomorphism
$$
\HH_G \otimes \E \cong \HH_G \otimes \F
$$
of $ G $-Hilbert $ B $-modules where $ \HH_G $ is viewed as a $ G $-Hilbert space using the left regular corepresentation, see \cite{Vergniouxthesis}.
Using the Kasparov stabilization theorem we deduce that
there is a $ G $-equivariant Hilbert $ B $-module isomorphism
$$
(\HH_G \otimes \E) \oplus (\HH_G \otimes \HH \otimes B) \cong \HH_G \otimes \HH \otimes B
$$
for every countably generated $ G $-Hilbert $ B $-module $ \E $. This result will be referred to as the equivariant stabilization theorem. \\
In the sequel we will frequently write $ KK^G $ instead of $ KK^S $ for $ S = C_0^\red(G) $ and call the defining cycles of this group $ G $-equivariant
Kasparov modules. It follows from Baaj-Skandalis duality that $ KK^G(A,B) $ can be represented by homotopy classes of $ G $-equivariant Kasparov
$ (\KH_G \otimes A) $-$ (\KH_G \otimes B) $-modules $ (\E, \phi, F) $ with $ G $-invariant operator $ F $. Taking Kasparov product with
the $ \KH_G \otimes B $-$ B $ imprimitivity bimodule $ (\HH_G \otimes B, \id, 0) $ we see that $ KK^G(A,B) $ can be represented by
homotopy classes of equivariant Kasparov $ (\KH_G \otimes A) $-$ B $ modules of the form $ (\HH_G \otimes \E, \phi, F) $ with invariant $ F $.
Using the equivariant stabilization theorem we can furthermore assume that $ (\HH_G \otimes \E)_\pm =
\HH_G \otimes \HH \otimes B $ is the standard $ G $-Hilbert $ B $-module. \\
From this point on we follow the arguments in \cite{Meyerkkg}.
Writing $ [A,B]_G $ for the set of equivariant homotopy classes of $ G $-equivariant $ * $-homomorphisms between $ G $-$ C^* $-algebras $ A $ and $ B $,
we arrive at the following description of the equivariant $ KK $-groups.
\begin{theorem} \label{cuntzpic}
Let $ G $ be a regular locally compact quantum group. Then there is a natural isomorphism
$$
KK^G(A,B) \cong [q(\KH_G \otimes A), \KH_G \otimes \KH \otimes B]_G
$$
for all separable $ G $-$ C^* $-algebras $ A $ and $ B $. We also have a natural isomophism
$$
KK^G(A,B) \cong [\KH_G \otimes \KH \otimes q(\KH_G \otimes \KH \otimes A), \KH_G \otimes \KH \otimes q(\KH_G \otimes \KH \otimes B)]_G
$$
under which the Kasparov product corresponds to the composition of homomorphisms.
\end{theorem}
Consider the category $ G\Alg $ of separable $ G $-$ C^* $-algebras for a regular quantum group $ G $. A functor $ F $ from $ G\Alg $
to an additive category $ \C $ is called a homotopy functor if $ F(f_0) = F(f_1) $ whenever $ f_0 $ and $ f_1 $ are $ G $-equivariantly
homotopic $ * $-homomorphisms. It is called stable if for all pairs of separable $ G $-Hilbert spaces $ \H_1, \H_2 $
the maps $ F(\KH(\H_j) \otimes A) \rightarrow F(\KH(\H_1 \oplus \H_2) \otimes A) $
induced by the canonical inclusions $ \H_j \rightarrow \H_1 \oplus \H_2 $ for $ j = 1,2 $
are isomorphisms. As in the group case, a homotopy functor
$ F $ is stable iff there exists a natural isomorphism
$ F(A) \cong F(\KH_G \otimes \KH \otimes A) $ for all $ A $.
Finally, $ F $ is called split exact if for every extension
\begin{equation*}
\xymatrix{
0 \ar@{->}[r] & K \ar@{->}[r] & E \ar@{->}[r] & Q \ar@{->}[r] & 0
}
\end{equation*}
of $ G $-$ C^* $-algebras that splits by an equivariant $ * $-homomorphism $ \sigma: Q \rightarrow E $
the induced sequence $ 0 \rightarrow F(K) \rightarrow F(E) \rightarrow F(Q) \rightarrow 0 $ in $ \C $ is split exact. \\
Equivariant $ KK $-theory can be viewed as an additive category $ KK^G $ with separable
$ G $-$ C^* $-algebras as objects and $ KK^G(A,B) $ as the set of morphisms between two objects $ A $ and $ B $.
Composition of morphisms is given by the Kasparov product.
There is a canonical functor $ \iota: G\Alg \rightarrow KK^G $ which is the identity on
objects and sends equivariant $ * $-homomorphisms to the corresponding $ KK $-elements. This functor is a split exact
stable homotopy functor. \\
As a consequence of theorem \ref{cuntzpic} we obtain the following universal property of $ KK^G $, see again \cite{Meyerkkg}.
We remark that a related assertion
is stated in \cite{Popescu}, however, some of the arguments in \cite{Popescu} are incorrect.
\begin{theorem} \label{kkuniversal}
Let $ G $ be a regular locally compact quantum group. The functor $ \iota: G \Alg \rightarrow KK^G $ is the universal split exact stable homotopy functor on the category $ G \Alg $.
More precisely, if $ F: G \Alg \rightarrow \C  $ is any split exact stable homotopy functor
with values in an additive category $ \C $ then there exists a unique functor $ f: KK^G \rightarrow \C $ such that
$ F = f \iota $.
\end{theorem}
Let us explain how $ KK^G $ becomes a triangulated category. We follow the discussion in \cite{MNtriangulated}, for the definition of
a triangulated category see \cite{Neeman}. Let $ \Sigma A $ denote the suspension $ C_0(\mathbb{R}) \otimes A $ of a $ G $-$ C^* $-algebra $ A $.
Here $ C_0(\mathbb{R}) $ is equipped with the trivial coaction. The corresponding functor $ \Sigma: KK^G \rightarrow KK^G $ determines the
translation automorphism.
If $ f: A \rightarrow B $ is a $ G $-equivariant $ * $-homomorphism then the mapping cone
$$
C_f = \{(a,b) \in A \times C_0((0,1], B) | b(1) = f(a) \}
$$
is a $ G $-$ C^* $-algebra in a natural way, and there is a canonical diagram
$$
\xymatrix{
\Sigma B  \;\; \ar@{->}[r] & C_f \ar@{->}[r] & A \ar@{->}[r]^f & B
}
$$
of $ G $-equivariant $ * $-homomorphisms. Diagrams of this form are called mapping cone triangles.
By definition, an exact triangle is a diagram $ \Sigma Q \rightarrow K \rightarrow E \rightarrow Q $ in $ KK^G $
which is isomorphic to a mapping cone triangle. \\
The proof of the following proposition is carried out in the same way as for locally compact groups \cite{MNtriangulated}.
\begin{prop}\label{kktriang}
Let $ G $ be a regular locally compact quantum group. Then the category $ KK^G $ together with the translation functor
and the exact triangles described above is triangulated.
\end{prop}
Several results about the equivariant $ KK $-groups for ordinary groups extend in a straightforward way
to the setting of quantum groups. As an example, let us state the Green-Julg theorem for compact quantum groups and its dual
version for discrete quantum groups.
If $ G $ is a locally compact quantum group and $ A $ is a $ C^* $-algebra we write $ \res^E_G(A) $ for the $ G $-$ C^* $-algebra $ A $
with the trivial coaction. A detailed proof of the following result is contained in \cite{Vergniouxthesis}.
\begin{theorem} \label{GJ}
Let $ G $ be a compact quantum group. Then there is a natural isomorphism
$$
KK^G(\res^E_G(A), B) \cong KK(A, C^*(G)^\cop \ltimes B)
$$
for all $ C^* $-algebras $ A $ and all $ G $-$ C^* $-algebras $ B $. \\
Dually, let $ G $ be a discrete quantum group. Then there is a natural isomorphism
$$
KK^G(A, \res^E_G(B)) \cong KK(C^*_\max(G)^\cop \ltimes_\max A, B)
$$
for all $ G $-$ C^* $-algebras $ A $ and all $ C^* $-algebras $ B $.
\end{theorem}
Let $ G $ be a strongly regular quantum group and let $ H \subset G $ be a regular closed quantum subgroup.
It is easy to check that restriction from $ G $ to $ H $ induces a triangulated functor $ \res^G_H: KK^G \rightarrow KK^H $.
This functor associates to a $ G $-$ C^* $-algebra $ A $ the
$ H $-$ C^* $-algebra $ \res^G_H(A) = A $ obtained by restricting the action.
Similarly, using the universal property of theorem \ref{kkuniversal} we obtain a triangulated functor $ \ind_H^G: KK^H \rightarrow KK^G $
which maps an $ H $-$ C^* $-algebra $ A $ to the induced $ G $-$ C^* $-algebra $ \ind_H^G(A) $. Note that the
compatibility of induction with stabilizations follows from Vaes' imprimitivity theorem stated above as theorem \ref{vaesimprimitivity}. \\
A closed quantum subgroup $ H \subset G $ is called cocompact if the $ C^* $-algebraic quantum homogeneous space $ C_0^\red(G/H) $ is a
unital $ C^* $-algebra. In this case we write $ C^\red(G/H) $ instead of $ C_0^\red(G/H) $.
Recall that a locally compact quantum group $ G $ is coamenable if the natural
map $ C^\max_0(G) \rightarrow C^\red_0(G) $ is an isomorphism.
Strong regularity is equivalent to regularity in this case.
\begin{prop} \label{Frobeniusrec}
Let $ H \subset G $ be a cocompact regular quantum subgroup of a strongly regular quantum group $ G $. If $ G $ is coamenable there is a
natural isomorphism
$$
KK^H(\res^G_H(A), B) \cong KK^G(A, \ind_H^G(B))
$$
for all $ G $-$ C^* $-algebras $ A $ and all $ H $-$ C^* $-algebras $ B $.
\end{prop}
\proof We describe the unit $ \eta $ and the counit $ \kappa $ of this adjunction.
For a $ G $-$ C^* $-algebra $ A $ let $ \eta_A: A \rightarrow \ind_H^G \res^G_H(A) \cong C^\red(G/H) \twisted_G A $
be the $ G $-equivariant $ * $-homomorphism obtained from the embedding of $ A $ in the braided tensor product.
Here we use theorem \ref{citensor} and the assumption that $ H \subset G $ is cocompact.
In order to define the counit $ \kappa $ recall that the induced $ C^* $-algebra $ \ind_H^G(B) $ of an $ H $-$ C^* $-algebra $ B $ is contained in
the $ C^\red_0(G) $-relative multiplier algebra of $ C^\red_0(G) \otimes B $.
We obtain an $ H $-equivariant $ * $-homomorphism $ \kappa_B: \res^G_H \ind_H^G(B) \rightarrow B $ as the restriction of $ \epsilon \otimes \id:
M(C^\red_0(G) \otimes B) \rightarrow M(B) $ where $ \epsilon: C^\red_0(G) \rightarrow \mathbb{C} $ is the counit. Here we use coamenability
of $ G $. \\
Let $ A $ be a $ G $-$ C^* $-algebra with coaction $ \alpha $ and let $ \res(\alpha): A \rightarrow M(C^\red_0(H) \otimes A) $
be the restriction of $ \alpha $ to $ H $. Using the relation
$ [(\I \otimes 1)\res(\alpha)(A)] = [\alpha(A) (\I \otimes 1)]  $ established in \cite{Vaesimprimitivity} we see that $ \kappa_{\res(A)} $ is given by
$ \epsilon \twisted \id: C^\red(G/H) \twisted_G A \rightarrow \mathbb{C} \twisted_G A \cong A  $.
Note that, although not being $ G $-equivariant, the map $ \epsilon $ is $ C^*_\red(G) $-colinear
and hence induces a $ * $-homomorphism between the braided tensor products as desired.
It follows that the composition
$$
\xymatrix{
\res^G_H(A) \;\; \ar@{->}[r]^{\!\!\!\!\!\!\!\!\!\!\!\! \res(\eta_A)} & \res^G_H \ind_H^G \res^G_H(A) \ar@{->}[r]^{\quad\quad\; \kappa_{\res(A)}} & \res^G_H(A)
}
$$
is the identity in $ KK^H(\res^G_H(A), \res^G_H(A)) $ for every $ G $-$ C^* $-algebra $ A $. \\
Identifying the isomorphism $ \ind^G_H \res_H^G \ind^G_H(B) \cong C^\red(G/H) \twisted_G \ind_H^G(B) $
and using the counit identity $ (\id \otimes \epsilon)\Delta = \id $ for $ C^\red_0(G) $ we see that
$$
\xymatrix{
\ind^G_H(B) \;\ar@{->}[r]^{\!\!\!\!\!\!\!\!\!\!\!\!\!\!\!\! \eta_{\ind(B)}} & \;\ind^G_H \res_H^G \ind^G_H(B)\; \ar@{->}[r]^{\quad\quad\;\, \ind(\kappa_B)} & \;\; \ind^G_H(B)
}
$$
is the identity in $ KK^G(\ind_H^G(B), \ind_H^G(B)) $ for every $ G $-$ C^* $-algebra $ B $. This yields the assertion. \qed \\
Based on the braided tensor product we introduce exterior products in equivariant Kasparov theory.
\begin{prop}\label{tauleft}
Let $ G $ be a regular locally compact quantum group, let $ A $ and $ B $ be $ G $-$ C^* $-algebras and let $ D $ be a
$ G $-$ \yd $-algebra. Then there exists a natural homomorphism
$$
\lambda_D: KK^G(A,B) \rightarrow KK^G(D \twisted_G A, D \twisted_G B)
$$
defining a triangulated functor $ \lambda_D: KK^G \rightarrow KK^G $. \\
If $ A $ and $ B $ are $ G $-$ \yd $-algebras then there is an analogous homomorphism
$$
\lambda_D: KK^{\DD(G)}(A,B) \rightarrow KK^{\DD(G)}(D \twisted_G A, D \twisted_G B)
$$
defining a triangulated functor $ \lambda_D: KK^{\DD(G)} \rightarrow KK^{\DD(G)} $.
\end{prop}
\proof We shall only discuss the first assertion, the case of $ G $-$ \yd $-algebras is treated
analogously. Taking the braided tensor product with $ D $ defines a split
exact homotopy functor from $ G \Alg $ to $ KK^G $. According to proposition \ref{twistedstability} this functor is stable.
Hence the existence of $ \lambda_D $ is a consequence of the universal property of $ KK^G $ established in theorem \ref{kkuniversal},
and the resulting functor is easily seen to be triangulated. \qed \\
The same arguments yield the following right-handed version of proposition \ref{tauleft}.
\begin{prop}\label{tauright}
Let $ G $ be a regular locally compact quantum group, let $ C $ and $ D $ be $ G $-$ \yd $-algebras and let
$ B $ be a $ G $-$ C^* $-algebra. Then there exists a natural homomorphism
$$
\rho_B: KK^{\DD(G)}(C,D) \rightarrow KK^G(C \twisted_G B, D \twisted_G B)
$$
defining a triangulated functor $ \rho_B: KK^{\DD(G)} \rightarrow KK^G $. \\
If $ B $ is a $ G $-$ \yd $-algebra we obtain a natural homomorphism
$$
\rho_B: KK^{\DD(G)}(C,D) \rightarrow KK^{\DD(G)}(C \twisted_G B, D \twisted_G B)
$$
defining a triangulated functor $ \rho_B: KK^{\DD(G)} \rightarrow KK^{\DD(G)} $.
\end{prop}
By construction, the class of a $ G $-equivariant $ * $-homomorphism $ f: A \rightarrow B $ is mapped to the class of
$ f \twisted \id $ under $ \lambda_D: KK^G \rightarrow KK^G $, and similar remarks apply to the other functors obtained above. \\
Of course one can also give direct definitions on the level of Kasparov modules for the constructions
in propositions \ref{tauleft} and \ref{tauright}. For instance, let $ (\E, \phi, F) $ be a 
$ G $-equivariant Kasparov $ A $-$ B $-module. Then $ D \twisted_G \E $ is a $ D \twisted_G B $-Hilbert module, 
and the map $ \phi: A \rightarrow \LH(\E) = M(\KH(\E)) $ induces a $ G $-equivariant $ * $-homomorphism $ \id \twisted_G \phi: 
D \twisted_G A \rightarrow M(D \twisted_G \KH(\E)) \cong \LH(D \twisted_G \E) $. Morever, we obtain $ \id \twisted_G F \in 
\LH(D \twisted_G \E) $ by applying the canonical map $ \LH(\E) \rightarrow \LH(D \twisted_G \E) $. It is readily checked that 
this yields a $ G $-equivariant Kasparov module. The construction is compatible with homotopies 
and induces $ \lambda_D: KK^G(A,B) \rightarrow KK^G(D \twisted_G A, D \twisted_G B) $. \\
Let $ A_1, B_1 $ and $ D $ be $ G $-$ \yd $ algebras and let $ A_2, B_2 $ be $ G $-$ C^* $-algebras. We define
the exterior Kasparov product
$$
KK^{\DD(G)}(A_1, B_1 \twisted_G D) \times KK^G(D \twisted_G A_2, B_2) \rightarrow KK^G(A_1 \twisted_G A_2, B_1 \twisted_G B_2)
$$
as the map which sends $ (x,y) $ to $ \rho_{A_2}(x) \circ \lambda_{B_1}(y) $. Here $ \circ $ denotes the
Kasparov composition product, and we use $ (B_1 \twisted_G D) \twisted_G A_2 \cong B_1 \twisted_G (D \twisted_G A_2) $. \\
If $ A_2, B_2 $ are $ G $-$ \yd $-algebras we obtain an exterior product
$$
KK^{\DD(G)}(A_1, B_1 \twisted_G D) \times KK^{\DD(G)}(D \twisted_G A_2, B_2) \rightarrow KK^{\DD(G)}(A_1 \twisted_G A_2, B_1 \twisted_G B_2)
$$
in the same way. \\
We summarize the main properties of the above exterior Kasparov products in analogy with the
ordinary exterior Kasparov product, see \cite{Blackadar}.
\begin{theorem} \label{exterior}
Let $ G $ be a regular locally compact quantum group. Moreover let $ A_1, B_1 $ and $ D $ be $ G $-$ \yd $ algebras and let $ A_2, B_2 $ be
$ G $-$ C^* $-algebras. The exterior Kasparov product
$$
KK^{\DD(G)}(A_1, B_1 \twisted_G D) \times KK^G(D \twisted_G A_2, B_2) \rightarrow KK^G(A_1 \twisted_G A_2, B_1 \twisted_G B_2)
$$
is associative and functorial in all possible senses. An analogous statement holds for the product
$$
KK^{\DD(G)}(A_1, B_1 \twisted_G D) \times KK^{\DD(G)}(D \twisted_G A_2, B_2) \rightarrow KK^{\DD(G)}(A_1 \twisted_G A_2, B_1 \twisted_G B_2)
$$
provided $ A_2, B_2 $ are $ G $-$ \yd $-algebras.
\end{theorem}
Recall that every $ G $-$ C^* $-algebra for a locally compact group $ G $ can be viewed as a $ G $-$ \yd $-algebra with the trivial
coaction of $ C^*_\red(G) $. In this case
our constructions reduce to the classical exterior product in equivariant $ KK $-theory. Still, even for classical
groups the products defined above are more general since we may consider $ G $-$ \yd $-algebras that are equipped
with a nontrivial coaction of the group $ C^* $-algebra.

\section{The quantum group $ SU_q(2) $} \label{secsuq2}

In this section we recall some definitions and constructions related to the compact quantum group $ SU_q(2) $.
For more information on the algebraic aspects of compact quantum groups we refer to \cite{KS}. \\
Let us fix a number $ q \in (0,1] $ and describe the $ C^* $-algebra of continuous functions on $ SU_q(2) $.
Since $ SU_q(2) $ is coamenable \cite{NagysuqN}, \cite{Banicafusion} there is no need to distinguish between the full and reduced $ C^* $-algebras.
By definition, $ C(SU_q(2)) $ is the universal $ C^* $-algebra generated by two elements $ \alpha $ and $ \gamma $ satisfying the relations
$$
\alpha \gamma = q \gamma \alpha, \quad \alpha \gamma^* = q \gamma^* \alpha, \quad \gamma \gamma^* = \gamma^* \gamma, \quad
\alpha^* \alpha + \gamma^* \gamma = 1, \quad \alpha \alpha^*  + q^2 \gamma \gamma^* = 1.
$$
The comultiplication $ \Delta: C(SU_q(2)) \rightarrow C(SU_q(2)) \otimes C(SU_q(2)) $ is given on
the generators by
$$
\Delta(\alpha) = \alpha \otimes \alpha - q \gamma^* \otimes \gamma, \qquad \Delta(\gamma) = \gamma \otimes \alpha + \alpha^* \otimes \gamma.
$$
From a conceptual point of view, it is useful to interpret these formulas in terms of the fundamental matrix
$$
u =
\begin{pmatrix}
\alpha & -q \gamma^* \\
\gamma & \alpha^*
\end{pmatrix}.
$$
In fact, the defining relations for $ C(SU_q(2)) $ are equivalent to saying that the fundamental matrix
is unitary, and the comultiplication of $ C(SU_q(2)) $ can be written in a concise
way as
$$
\Delta
\begin{pmatrix}
\alpha & -q \gamma^* \\
\gamma & \alpha^*
\end{pmatrix}
=
\begin{pmatrix}
\alpha & -q \gamma^* \\
\gamma & \alpha^*
\end{pmatrix}
\otimes
\begin{pmatrix}
\alpha & -q \gamma^* \\
\gamma & \alpha^*
\end{pmatrix}.
$$
We will also work with the dense
$ * $-subalgebra $ \mathbb{C}[SU_q(2)] \subset C(SU_q(2)) $ generated by $ \alpha $ and $ \gamma $.
Together with the counit $ \epsilon: \mathbb{C}[SU_q(2)] \rightarrow \mathbb{C} $ and the
antipode $ S: \mathbb{C}[SU_q(2)] \rightarrow \mathbb{C}[SU_q(2)] $ determined by
$$
\epsilon\begin{pmatrix}
\alpha & -q \gamma^* \\
\gamma & \alpha^*
\end{pmatrix}
=
\begin{pmatrix}
1 & 0 \\
0 & 1
\end{pmatrix},
\qquad
S\begin{pmatrix}
\alpha & -q \gamma^* \\
\gamma & \alpha^*
\end{pmatrix}
=
\begin{pmatrix}
\alpha^* & \gamma^* \\
-q \gamma & \alpha
\end{pmatrix}
$$
the algebra $ \mathbb{C}[SU_q(2)] $ becomes a Hopf-$ * $-algebra.
We use the Sweedler notation $ \Delta(x) = x_{(1)} \otimes x_{(2)} $ for the
comultiplication and write
$$
f \hit x = x_{(1)} f(x_{(2)}), \qquad x \hitby f = f(x_{(1)}) x_{(2)}
$$
for elements $ x \in \mathbb{C}[SU_q(2)] $ and linear functionals $ f: \mathbb{C}[SU_q(2)] \rightarrow \mathbb{C} $. \\
The antipode is an algebra antihomomorphism satisfying
$ S(S(x^*)^*) = x $ for all $ x \in \mathbb{C}[SU_q(2)] $, in particular the map $ S $ is invertible.
The inverse of $ S $ can be written as
$$
S^{-1}(x) = \delta \hit S(x) \hitby \delta^{-1}
$$
where $ \delta: \mathbb{C}[SU_q(2)] \rightarrow \mathbb{C} $ is the modular character determined by
$$
\delta\begin{pmatrix}
\alpha & -q \gamma^* \\
\gamma & \alpha^*
\end{pmatrix}
=
\begin{pmatrix}
q^{-1} & 0 \\
0 & q
\end{pmatrix}.
$$
Apart from its role in connection with the antipode, the character $ \delta $ describes the modular properties of the Haar state $ \phi $ of $ C(SU_q(2)) $ in the sense that
$$
\phi(xy) = \phi(y (\delta \hit x \hitby \delta))
$$
for all $ x,y \in \mathbb{C}[SU_q(2)] $.
The Hilbert space $ \HH_{SU_q(2)} $ associated to $ SU_q(2) $ is the GNS-construction of $ \phi $ and will be denoted by $ L^2(SU_q(2)) $
in the sequel. \\
The irreducible corepresentations $ V_l $ of $ C(SU_q(2)) $ are parametrized by $ l \in \frac{1}{2} \mathbb{N} $,
and the dimension of $ V_l $ is $ 2l + 1 $ as for the classical group $ SU(2) $. According to the Peter-Weyl theorem, the Hilbert space $ L^2(SU_q(2)) $
has an orthonormal basis $ e^{(l)}_{ij} $ with $ l \in \frac{1}{2} \mathbb{N} $ and
$ i,j \in \{-l, -l + 1, \dots, l \} $ corresponding to the decomposition of the regular corepresentation.
In this picture, the GNS-representation of $ C(SU_q(2)) $ is given by
\begin{align*}
\alpha \, e_{ij}^{(l)} &= a_+(l,i,j) \, e_{i - \frac{1}{2}, j - \frac{1}{2}}^{\left(l + \frac{1}{2} \right)}
+ a_- (l,i,j) \, e_{i - \frac{1}{2} , j - \frac{1}{2}}^{\left(l - \frac{1}{2} \right)} \\
\gamma \, e_{ij}^{(l)} &= c_+ (l,i,j) \, e_{i + \frac{1}{2}, j - \frac{1}{2}}^{\left(l + \frac{1}{2} \right)} +
c_- (l,i,j) \, e_{i + \frac{1}{2} , j - \frac{1}{2}}^{\left(l - \frac{1}{2} \right)}
\end{align*}
where the explicit form of $ a_\pm $ and $ c_\pm $ for $ q \in (0,1) $ is
\begin{align*}
a_+ (l,i,j) &= q^{2l + i + j + 1} \, \frac{(1 - q^{2l - 2j + 2})^{1/2} (1 - q^{2l - 2i + 2})^{1/2}}{ (1 - q^{4l + 2})^{1/2} (1 - q^{4l + 4})^{1/2}} \\
a_- (l,i,j) &= \frac{(1 - q^{2l + 2j})^{1/2} (1 - q^{2l + 2i})^{1/2}}{(1 - q^{4l})^{1/2} (1 - q^{4l+2})^{1/2}}
\end{align*}
and
\begin{align*}
c_+ (l,i,j) &= -q^{l + j} \, \frac{(1 - q^{2l - 2j + 2})^{1/2}
(1 - q^{2l + 2i + 2})^{1/2}}{(1 - q^{4l + 2})^{1/2} (1 - q^{4l + 4})^{1/2}} \\
c_- (l,i,j) &= q^{l + i} \, \frac{(1 - q^{2l + 2j})^{1/2}
(1 - q^{2l - 2i})^{1/2}}{(1 - q^{4l})^{1/2} (1 - q^{4l + 2})^{1/2}}.
\end{align*}
In the above formulas the vectors $ e^{(l)}_{ij} $ are declared to be zero if one of the indices $ i, j $ is not
contained in $ \{-l, -l + 1, \dots, l \} $. \\
We will frequently use the fact that the classical torus $ T = S^1 $ is a closed quantum subgroup of $ SU_q(2) $.
The inclusion $ T \subset SU_q(2) $ is determined by the $ * $-homomorphism
$ \pi: \mathbb{C}[SU_q(2)] \rightarrow \mathbb{C}[T] = \mathbb{C}[z,z^{-1}] $ given by
$$
\pi
\begin{pmatrix}
\alpha & -q \gamma^* \\
\gamma & \alpha^*
\end{pmatrix}
=
\begin{pmatrix}
z & 0 \\
0 & z^{-1}
\end{pmatrix}.
$$
By definition, the standard Podle\'s sphere $ C(SU_q(2)/T)  $ is the corresponding homogeneous space.
In the algebraic setting, the Podle\'s sphere is described by the dense
$ * $-subalgebra $ \mathbb{C}[SU_q(2)/T] \subset C(SU_q(2)/T) $ of coinvariants in $ \mathbb{C}[SU_q(2)] $
with respect to the right coaction $ (\id \otimes \pi)\Delta $ of $ \mathbb{C}[T] $. \\
If $ V $ is a finite dimensional left $ \mathbb{C}[T] $-comodule, or equivalently a finite dimensional
representation of $ T $, then the cotensor product
\begin{equation*}
\Gamma(SU_q(2) \times_T V) = \mathbb{C}[SU_q(2)] \Box_{\mathbb{C}[T]} V \subset \mathbb{C}[SU_q(2)] \otimes V
\end{equation*}
is a noncommutative analogue of the space of sections of the homogeneous vector bundle $ SU(2) \times_T V $ over $ SU(2)/T $.
Clearly $ \Gamma(SU_q(2) \times_T V) $ is a $ \mathbb{C}[SU_q(2)/T] $-bimodule in a natural way. In accordance with the Serre-Swan
theorem, the space of sections $ \Gamma(SU_q(2) \times_T V) $ is finitely generated and projective
both as a left and right $ \mathbb{C}[SU_q(2)/T]$-module. This follows from the fact that
$ \mathbb{C}[SU_q(2)/T] \subset \mathbb{C}[SU_q(2)] $ is a faithfully flat Hopf-Galois extension, see \cite{MS}, \cite{S}.
If $ V = \mathbb{C}_k $ is the irreducible representation of $ T $ of weight $ k \in \mathbb{Z} $ we write
$ L^2(SU_q(2) \times_T \mathbb{C}_k) $ for the $ SU_q(2) $-Hilbert space obtained by taking the closure of $ \Gamma(SU_q(2) \times_T \mathbb{C}_k) $
inside $ L^2(SU_q(2)) $. We also note the Frobenius reciprocity isomorphism
$$
\Hom_T(\res^{SU_q(2)}_T(V), \mathbb{C}_k) \cong \Hom_{SU_q(2)}(V, L^2(SU_q(2) \times_T \mathbb{C}_k))
$$
for all finite dimensional corepresentations $ V $ of $ C(SU_q(2)) $.

\section{Equivariant Poincar\'e duality for the Podle\'s sphere} \label{secfredpodles}

Poincar\'e duality in Kasparov theory plays an important r\^ole in noncommutative geometry,
for instance in connection with the Dirac-dual Dirac method for proving the Novikov conjecture \cite{Kasparov2}.
In this section we extend this concept to the setting of quantum group actions
and show that the standard Podle\'s sphere is equivariantly Poincar\'e dual to itself. \\
Let us begin with the following terminology, generalizing the definition given by Connes in \cite{Connes2}. Recall that we write $ \DD(G) $ for the
Drinfeld double of a locally compact quantum group $ G $.
\begin{definition} \label{PDdef}
Let $ G $ be a regular locally compact quantum group. Two $ G $-$ \yd $-algebras $ P $ and $ Q $ are called $ G $-equivariantly Poincar\'e dual to each other
if there exists a natural isomorphism
$$
KK^{\DD(G)}_*(P \twisted_G A, B) \cong KK^{\DD(G)}_*(A, Q \twisted_G B)
$$
for all $ G $-$ \yd $-algebras $ A $ and $ B $.
\end{definition}
Using the notation introduced in proposition \ref{tauleft} we may rephrase this by saying that the $ G $-$ \yd $-algebras $ P $ and $ Q $ are $ G $-equivariantly Poincar\'e dual to each other iff $ \lambda_P $ and $ \lambda_Q $ are adjoint functors. In particular, the unit and counit of the adjunction
determine elements
$$
\alpha \in KK^{\DD(G)}_*(P \twisted_G Q, \mathbb{C}), \qquad \beta \in KK^{\DD(G)}_*(\mathbb{C}, Q \twisted_G P)
$$
if $ P $ and $ Q $ are Poincar\'e dual. In this case one also has a duality on the level of $ G $-equivariant Kasparov theory in the sense
that there is a natural isomorphism
$$
KK^G_*(P \twisted_G A, B) \cong KK^G_*(A, Q \twisted_G B)
$$
for all $ G $-$ C^* $-algebras $ A $ and $ B $. \\
In the sequel we restrict attention to $ G_q = SU_q(2) $. Our aim is to show that the standard Podle\'s sphere is $ SU_q(2) $-equivariantly Poincar\'e dual to
itself in the sense of definition \ref{PDdef}. As a first ingredient we need the $ K $-homology class of the Dirac operator on $ G_q/T $
for $ q \in (0,1) $. We review briefly the construction in \cite{DSPodles}, however, instead of working with the action of the quantized universal
enveloping algebra
we consider the corresponding coaction of $ C(G_q) $. Using the notation from section \ref{secsuq2}, the underlying
graded $ G_q $-Hilbert space $ \H = \H_+ \oplus \H_- $ of the spectral triple is
given by
$$
\H_\pm = L^2(G_q \times_T \mathbb{C}_{\pm 1})
$$
with its natural coaction of $ C(G_q) $. The covariant representation $ \phi = \phi_+ \oplus \phi_- $ of the $ C(G_q/T) $
is given by left multiplication. Finally, the Dirac operator $ D $ on $ \H $ is the odd operator
$$
D =
\begin{pmatrix}
0 & D^- \\
D^+ & 0
\end{pmatrix}
$$
where
$$
D^\pm|l,m \ket_\pm = [l + 1/2]_q\; |l,m \ket_\mp
$$
and $ |l,m \ket_\pm $ are the standard basis vectors in $ V_l \subset \H_\pm $ and
$$
[a]_q = \frac{q^a - q^{-a}}{q - q^{-1}}
$$
for a nonzero number $ a \in \mathbb{C} $.
Note that $ \H_+ $ and $ \H_- $ are isomorphic corepresentations of
$ C(G_q) $ according to Frobenius reciprocity. It follows that the phase $ F $ of $ D $ can be written as
$$
F =
\begin{pmatrix}
0 & 1 \\
1 & 0
\end{pmatrix},
$$
and the triple $ (\H, \phi, F) $ is a $ G_q $-equivariant Fredholm module. In this way $ D $
determines an element in $ KK^{G_q}_0(C(G_q/T), \mathbb{C}) $. \\
According to proposition \ref{indyd} the $ C^* $-algebra  $ C(G_q/T) = \ind_T^{G_q}(\mathbb{C}) $ is a $ G_q $-$ \yd $-algebra.
For our purposes the following fact is important.
\begin{prop} \label{DiracYD}
The Dirac operator on the standard Podle\'s sphere defines an element in $ KK^{\DD(G_q)}(C(G_q/T) \twisted_{G_q} C(G_q/T), \mathbb{C}) $ in a natural way.
\end{prop}
\proof With the notation as above, we consider the operator $ F $ on the Hilbert space $ \H = \H_+ \oplus \H_- $. Using
$ \ind_T^{G_q }\res^{G_q}_T(C(G_q/T)) \cong C(G_q/T) \twisted_{G_q} C(G_q/T) $ we
obtain a graded $ G_q $-equivariant $ * $-homomorphism $ \psi: C(G_q/T) \twisted_{G_q} C(G_q/T) \rightarrow \LH(\H) $
by applying the induction functor to the counit $ \epsilon: C(G_q/T) \rightarrow \mathbb{C} $ and composing the resulting map
with the natural representation of $ C(G_q/T) $ on $ \H $. On both copies
of $ C(G_q/T) $ the map $ \psi $ is given by the homomorphism $ \phi $ from above.
In particular, the commutators of $ F $ with elements from
$ C(G_q/T) \twisted_{G_q} C(G_q/T) $ are compact. \\
The coaction $ \lambda: \H \rightarrow M(C^*(G_q) \otimes \H) $ which turns $ \H $ into a $ \DD(G_q) $-Hilbert space
is obtained from the action of $ \mathbb{C}[G_q] $ on $ \Gamma(G_q \times_T \mathbb{C}_\pm) $ given by
$$
f \cdot h = f_{(1)} h \delta \hit S(f_{(2)})
$$
where $ \delta $ is the modular character. The homomorphism $ \psi $ is $ C^*(G_q) $-colinear with respect to this
coaction, and in order to show
$$
(C^*(G_q) \otimes 1)(1 \otimes F - \ad_\lambda(F)) \subset C^*(G_q) \otimes \KH(\H)
$$
it suffices to check that $ F $ commutes with the above action of $ \mathbb{C}[G_q] $ up to compact operators. This in turn is
a lengthy but straightforward calculation based on the explicit formulas for the GNS-representation of $ C(G_q) $
in section \ref{secsuq2}. It follows that $ (\H, \psi, F) $ is a $ \DD(G_q) $-equivariant Kasparov module as desired. \qed \\
Note that in the construction of the Dirac cycle in proposition \ref{DiracYD} we use two identical representations of $ C(G_q/T) $
as in the case of a classical spin manifold. The difference to the classical situation lies in the replacement of the ordinary tensor
product with the braided tensor product. \\
Let us formally write $ \E_k = G_q \times_T \mathbb{C}_k $ for the induced vector bundle associated to the representation
of weight $ k $, and denote by $ C(\E_k) $ the closure of $ \Gamma(\E_k) $ inside $ C(G_q) $.
The space $ C(\E_k) $ is a $ G_q $-equivariant Hilbert $ C(G_q/T) $-module with the coaction induced by comultiplication,
and the coaction $ \lambda: C(\E_k) \rightarrow M(C^*(G_q) \otimes C(\E_k)) $
given by $ \lambda(f) = \hat{W}^*(1 \otimes f) \hat{W} $ turns it into a $ \DD(G_q) $-equivariant Hilbert module.
Left multiplication yields a $ \DD(G_q) $-equivariant $ * $-homomorphism $ \mu: C(G_q/T) \rightarrow \KH(C(\E_k)) $.
Hence $ (C(\E_k), \mu, 0) $ defines a class $ [[\E_k]] $ in
$ KK^{\DD(G_q)}_0(C(G_q/T), C(G_q/T)) $. \\
Next observe that the unit homomorphism $ u: \mathbb{C} \rightarrow C(G_q/T) $ induces an element $ [u] \in KK^{\DD(G_q)}_0(\mathbb{C}, C(G_q/T)) $.
We obtain a class $ [\E_k] $ in $ KK^{\DD(G_q)}_0(\mathbb{C}, C(G_q/T)) $ by restricting $ [[\E_k]] $ along $ u $,
or equivalently, by taking the product
$$
[\E_k] = [u] \circ [[\E_k]].
$$
Under the forgetful map from $ KK^{\DD(G_q)} $ to $ KK^{G_q} $, this class is
mapped to the $ K $-theory class in $ KK^{G_q}(\mathbb{C}, C(G_q/T)) $ corresponding
to $ \mathbb{C}_k $ in $ R(T) $ under Frobenius reciprocity. \\
In addition we define elements $ [D \otimes \E_k] \in KK^{G_q}_0(C(G_q/T), \mathbb{C}) $ by
$$
[D \otimes \E_k] = [[\E_k]] \circ [D]
$$
where $ [D] \in KK^{G_q}_0(C(G_q/T), \mathbb{C}) $ is the class of the
Dirac operator. We remark that these elements correspond to twisted Dirac operators on $ G_q/T $ as studied by Sitarz in \cite{Sitarztwisted}. \\
Let us determine the equivariant indices of these twisted Dirac operators.
\begin{prop} \label{indexhomogenous}
Consider the classes $ [\E_k] \in KK^{G_q}(\mathbb{C}, C(G_q/T)) $ and $ [D \otimes \E_l] \in KK^{G_q}(C(G_q/T), \mathbb{C}) $
introduced above. The Kasparov product $ [\E_k] \circ [D \otimes \E_l] $ in $ KK^{G_q}_0(\mathbb{C}, \mathbb{C}) = R(G_q) $
is given by
$$
[\E_k] \circ [D \otimes \E_l] =
\begin{cases}
-[V_{(k + l - 1)/2}] & \quad \text{for}\; k + l > 0 \\
0 & \quad \text{for}\; k + l = 0 \\
[V_{-(k + l + 1)/2}] & \quad \text{for}\; k + l < 0
\end{cases}
$$
for $ k, l \in \mathbb{Z} $.
\end{prop}
\proof This is analogous to calculating the index of a homogenous differential operator \cite{Bott}.
Since we have $ [[\E_m]] \circ [[\E_n]] = [[\E_{m + n}]] $ for all $ m,n \in \mathbb{Z} $ it suffices to consider the case
$ k = 0 $. The product $ [\E_0] \circ [D \otimes \E_l] $ is given by the equivariant index of the $ G_q $-equivariant Fredholm operator
representing $ [D \otimes \E_l] $. This operator can be viewed
as an odd operator on $ L^2(\E_{l + 1}) \oplus L^2(\E_{l - 1}) $. By equivariance, the claim follows
from Frobenius reciprocity; we only have to subtract the classes of $ L^2(\E_{l + 1}) $ and $ L^2(\E_{l - 1}) $
in the formal representation ring of $ G_q $. \qed \\
We note that for the above computation there is no need to pass to cyclic cohomology or twisted cyclic cohomology. \\
In order to proceed we need a generalization of the Drinfeld double.
The relative Drinfeld double $ \DD(T, \hat{G}_q) $ is defined as the double crossed product \cite{BV} of $ C(T) $ and $ C^*(G_q) $
using the matching $ m(x) = \mathcal{Z} x \mathcal{Z}^* $ where $ \mathcal{Z} = (\pi \otimes \id)(W_{G_q}) $ and
$ \pi: C(G_q) \rightarrow C(T) $ is the quotient map. That is, we have
$ C^\red_0(\DD(T, \hat{G}_q)) = C(T) \otimes C^*(G_q) $ with the comultiplication
$$
\Delta_{\DD(T, \hat{G}_q)} = (\id \otimes \sigma \otimes \id)(\id \otimes m \otimes \id)(\Delta \otimes \hat{\Delta}).
$$
The relative Drinfeld double $ \DD(T, \hat{G}_q) $ is a cocompact closed quantum subgroup of $ \DD(G_q) $, and the quantum homogeneous
space $ C^\red(\DD(G_q)/\DD(T, \hat{G}_q)) $ is isomorphic to $ C(G_q/T) $. Under this identification, the natural $ \DD(G_q) $-algebra structure on the
homogeneous space corresponds to the $ \yd $-algebra structure on the induced algebra $ C(G_q/T) = \ind_T^{G_q}(\mathbb{C}) $
obtained from proposition \ref{indyd}. \\
Every continuous coaction of $ C(T) $ on a $ C^* $-algebra $ B $ restricts to a continuous coaction of $ C^\red_0(\DD_q) = C^\red_0(\DD(T, \hat{G}_q)) $
in a natural way, and we write $ \res^T_{\DD_q}(B) $ for the resulting $ \DD(T, \hat{G}_q) $-$ C^* $-algebra. Indeed, since $ C(T) $ is commutative,
the canonical $ * $-homomorphism $ C(T) \rightarrow M(C^\red_0(\DD(T, \hat{G}_q)) $ is compatible with the comultiplications.
The following result is a variant of the dual Green-Julg theorem, see theorem \ref{GJ}.
\begin{lemma} \label{DGJextended}
Let $ \DD_q = \DD(T, \hat{G}_q) $ be the relative Drinfeld double of $ G_q $. Then there is a natural isomorphism
$$
KK^{\DD_q}(A, \res^T_{\DD_q}(B)) \cong KK^T(C(G_q)^\cop \ltimes A, B)
$$
for all $ \DD_q $-$ C^* $-algebras $ A $ and all $ T $-$ C^* $-algebras $ B $.
\end{lemma}
\proof If $ A $ is a $ \DD_q $-$ C^* $-algebra then the crossed product $ C(G_q)^\cop \ltimes A $ becomes a $ T $-$ C^* $-algebra using
the adjoint action on $ C(G_q) $ and the restriction of the given coaction on $ A $. The natural map
$ \iota_A: A \rightarrow C(G_q)^\cop \ltimes A $ is $ T $-equivariant and $ C^*(G_q) $-colinear
with respect to the coaction on the crossed product induced by the corepresentation $ \hat{W}_G $. \\
Assume that $ (\E, \phi, F) $ is a $ \DD_q $-equivariant Kasparov $ A $-$ \res^T_{\DD_q}(B) $-module which is essential in the sense that
the $ * $-homomorphism $ \phi: A \rightarrow \LH(\E) $ is nondegenerate. The coaction of $ C^\red_0(\DD_q) $ on $ \E $ is determined by a
coaction of $ C(T) $ and a unitary corepresentation of $ C^*(G_q) $. Together with $ \phi $, this corepresentation
corresponds to a nondegenerate $ * $-homomorphism $ \psi: C(G_q)^\cop \ltimes A \rightarrow \LH(\E) $ which
yields a $ T $-equivariant Kasparov $ C(G_q)^\cop \ltimes A $-$ B $-module $ (\E, \psi, F) $.
Conversely, assume that $ (\E, \psi, F) $ is an essential $ T $-equivariant Kasparov $ C(G_q)^\cop \ltimes A $-$ B $-module.
Then $ \psi $ is determined by a covariant pair consisting of a nondegenerate $ * $-homomorphism $ \phi: A \rightarrow \LH(\E) $
and a unitary corepresentation of $ C^*(G_q) $ on $ \E $.
In combination with the given $ C(T) $-coaction, this corepresentation determines a coaction of $ C^\red_0(\DD_q) $ on $ \E $ such that
$ (\E, \phi, F) $ is a $ \DD_q $-equivariant Kasparov module. The assertion follows easily from these observations. \qed \\
Before we proceed we need some further facts about the structure of $ q $-deformations.
Note that $ C(G_q) $ can be viewed as a $ T \times T $-$ C^* $-algebra
with the action given by left and right translations. The $ C^* $-algebras $ C(G_q) $
assemble into a $ T \times T $-equivariant continuous field $ {\bf G} = (C(G_q))_{q \in (0,1]} $
of $ C^* $-algebras, compare \cite{Blanchard}, \cite{Nagydeformation}. In particular, the algebra $ C_0({\bf G}) $ of $ C_0 $-sections of the field
is a $ T \times T $-$ C^* $-algebra in a natural way.
We can also associate equivariant continuous fields to certain braided tensor products. For instance,
the braided tensor products $ C(G_q) \twisted_{G_q} C(G_q) $ yield
a continuous field of $ C^* $-algebras over $ (0,1] $ whose section algebra we denote by
$ C_0({\bf G}) \twisted_{\bf G} C_0({\bf G}) $. This is
easily seen using that $ C(G_q) \twisted_{G_q} C(G_q) \cong C(G_q) \otimes C(G_q) $
as $ C^* $-algebras and the fact that $ C(G_q) $ is nuclear for all $ q \in (0,1] $. A similar argument
works for the quantum flag manifolds $ C(G_q/T) $ instead of $ C(G_q) $. \\
As a consequence of lemma \ref{DGJextended} we obtain in particular that the Dirac operator on $ G_q/T $ determines an element in
$ KK^T(C(G_q/T) \twisted_{G_q} C(G_q/T) \twisted_{G_q} C(G_q), \mathbb{C}) $ since we have
$$
C(G_q)^\cop \ltimes A \cong A \twisted_{G_q} C(G_q)
$$
for every $ G_q $-$ \yd $-algebra $ A $ by definition of the braided tensor product. \\
In fact, these elements depend in a continuous way on the deformation parameter. More precisely, if we fix $ q \in (0,1] $
then the proof of proposition \ref{DiracYD} shows that the Dirac operators
on $ G_t/T $ for different values of $ t \in [q,1] $ yield an element
$$
[{\bf D}] \in KK^T(C({\bf G}/T) \twisted_{\bf G} C({\bf G}/T) \twisted_{\bf G} C({\bf G}), C[q,1])
$$
where $ C({\bf G}/T) \twisted_{\bf G} C({\bf G}/T) \twisted_{\bf G} C({\bf G}) $ denotes the algebra of sections of the continuous field
over $ [q,1] $ with fibers $ C(G_t/T) \twisted_{G_t} C(G_t/T) \twisted_{G_t} C(G_t) $. \\
Similarly, writing $ C({\bf G}/T) $ for the algebra of sections of the continuous field over $ [q,1] $
given by the Podle\'s spheres, the induced vector bundle $ \E_k $ determines a class in
$ KK^T(C({\bf G}/T), C({\bf G}/T)) $. Composition of this class with the canonical homomorphism $ C[q,1] \rightarrow C({\bf G}/T) $
yields an element in $ KK^T(C[q,1], C({\bf G}/T)) $. \\
After these preparations we prove the following main result.
\begin{theorem} \label{PD}
The Podle\'s sphere $ C(G_q/T) $ is $ G_q $-equivariantly Poincar\'e dual to itself. That is, there is a natural
isomorphism
$$
KK^{\DD(G_q)}_*(C(G_q/T) \twisted_{G_q} A, B) \cong KK^{\DD(G_q)}_*(A, C(G_q/T) \twisted_{G_q} B)
$$
for all $ G_q $-$ \yd $-algebras $ A $ and $ B $.
\end{theorem}
\proof According to proposition \ref{DiracYD} the Dirac operator on $ G_q/T $ yields an element $ [D_q]
\in KK^{\DD(G_q)}_0(C(G_q/T) \twisted_{G_q} C(G_q/T), \mathbb{C}) $.
Let us define a dual element $ \eta_q $ in $ KK^{\DD(G_q)}_0(\mathbb{C}, C(G_q/T) \twisted_{G_q} C(G_q/T)) $ by
$ \eta_q = [\E_{-1}] \twisted [\E_0] - [\E_0] \twisted [\E_1] $ where we write $ \twisted $ for
the exterior product obtained in theorem \ref{exterior}. \\
In order to show that $ \eta_q $ and $ [D_q] $ are the unit and counit of the desired adjunction we have to study the
endomorphisms $ (\id \twisted \eta_q) \circ ([D_q] \twisted \id) $ and $ (\eta_q \twisted \id) \circ (\id \twisted [D_q]) $
of $ C(G_q/T) $ in $ KK^{\DD(G_q)} $. \\
First we consider the classical case $ q = 1 $. Since all $ C^*_\red(G_1) $-coactions in the construction of $ [D_1] $ and $ \eta_1 $ are trivial
it suffices to work with the above morphisms at the level of $ KK^{G_1} $.
Due to proposition \ref{Frobeniusrec} the counit $ \epsilon: C(G_q/T) \rightarrow \mathbb{C} $ induces an
isomorphism $ KK^{G_1}_*(C(G_1/T), C(G_1/T)) \cong KK^T_*(C(G_1/T), \mathbb{C}) $.
Hence, according to the universal coefficient theorem for $ T $-equivariant $ KK $-theory \cite{RS}, in order to identify
$ (\id \twisted \eta_1) \circ ([D_1] \twisted \id) $ we only have to compute the action of
$ (\id \twisted \eta_1) \circ ([D_1] \twisted \id) \circ \epsilon $ on $ K^T_*(C(G_1/T)) $. Using proposition \ref{indexhomogenous} we
obtain
\begin{align*}
[\E_0] \circ (\id \twisted \eta_1) &\circ ([D_1] \twisted \id) \circ \epsilon =
[\E_0] \circ (\id \twisted \eta_1) \circ (\id \twisted \id \twisted \epsilon) \circ [D_1]\\
&= (([\E_0] \twisted [\E_{-1}]) \circ [D_1])\twisted [\mathbb{C}_0] - ([\E_0] \twisted [\E_0]) \circ [D_1]) \twisted [\mathbb{C}_1] = [\mathbb{C}_0]
\end{align*}
in $ K^T_0(\mathbb{C}) = R(T) $. Similarly
one checks $ [\E_1] \circ (\id \twisted \eta_1) \circ ([D_1] \twisted \id) \circ \epsilon = [\mathbb{C}_1] $.
This implies $ (\id \twisted \eta_1) \circ ([D_1] \twisted \id) = \id $
since $ [\E_0] $ and $ [\E_1] $ generate $ K^T_*(G_1/T) \cong R(T) \otimes_{R(G_1)} R(T) $ due to McLeod's theorem \cite{McLeod}.
In a similar way one shows $ (\eta_1 \twisted \id) \circ (\id \twisted [D_1]) = \id $.
As already indicated above, we conclude that these identities hold at the level of $ KK^{\DD(G_1)} $ as well. \\
For general $ q \in (0,1] $ we observe that the Drinfeld double $ \DD(G_q) $ is coamenable
and recall that $ \DD(T, \hat{G}_q) \subset \DD(G_q) $ is a cocompact quantum subgroup. According to proposition \ref{Frobeniusrec}
this implies
\begin{align*}
KK^{\DD(G_q)}_*(C(G_q/T), C(G_q/T)) \cong KK^{\DD(T, \hat{G}_q)}_*(C(G_q/T), \mathbb{C})
\end{align*}
since $ C(G_q/T) \cong \ind_{\DD(T, \hat{G}_q)}^{\DD(G_q)}(\mathbb{C}) $. Moreover,
due to lemma \ref{DGJextended} we have
\begin{align*}
KK^{\DD(T, \hat{G}_q)}_*(C(G_q/T), \mathbb{C}) \cong KK^T_*(C(G_q/T) \twisted_{G_q} C(G_q), \mathbb{C})
\end{align*}
using $ C(G_q)^\cop \ltimes C(G_q/T) \cong C(G_q/T) \twisted_{G_q} C(G_q) $. Recall that $ T $ acts by conjugation
on the copy of $ C(G_q) $. \\
The element in $ KK^T_*(C(G_q/T) \twisted_{G_q} C(G_q), \mathbb{C}) $ corresponding to
$ (\id \twisted \eta_q) \circ ([D_q] \twisted \id) $ is given by
$$
\delta_q = (\id \twisted \eta_q \twisted \id) \circ ([D_q] \twisted \id \twisted \id) \circ (\epsilon \twisted \id) \circ \epsilon.
$$
We observe that the individual elements in this composition assemble into $ KK^T $-classes
for the corresponding continuous fields over $ [q,1] $. \\
Let us denote by $ c_q \in E^T_0(C(G_1/T) \twisted_{G_1} C(G_1), C(G_q/T) \twisted_{G_q} C(G_q)) $ the $ E $-theoretic
comparison element for the field $ C({\bf G}/T) \twisted_{{\bf G}} C({\bf G}) $ over $ [q,1] $.
Using again the universal coefficient theorem for $ T $-equivariant $ KK $-theory we obtain a commutative diagram
$$
\xymatrix{
C(G_1/T) \twisted_{G_1} C(G_1) \ar@{->}[r]^{c_q} \ar@{->}[d]^{\delta_1} & C(G_q/T) \twisted_{G_q} C(G_q)
\ar@{->}[d]^{\delta_q} \\
\mathbb{C} \ar@{->}[r]^\id & \mathbb{C}
     }
$$
in $ KK^T $ where $ c_q $ is an isomorphism. Moreover, due to our previous considerations in the case $ q = 1 $
we have $ \delta_1 = (\epsilon \twisted \id) \circ \epsilon $.  This implies $ \delta_q = (\epsilon \twisted \id) \circ \epsilon $ and
hence $ (\id \twisted \eta_q) \circ ([D_q] \twisted \id) = \id $ in $ KK^{\DD(G_q)} $. In a similar way one
obtains $ (\eta_q \twisted \id) \circ (\id \twisted [D_q]) = \id $ in $ KK^{\DD(G_q)} $. According to the
characterization of adjoint functors in terms of unit and counit this yields the assertion. \qed \\
As a corollary we determine the equivariant $ K $-homology of the Podle\'s sphere.
\begin{cor}
For the standard Podle\'s sphere $ C(G_q/T) $ we have
$$
KK^{G_q}_0(C(G_q/T), \mathbb{C}) \cong R(G_q) \oplus R(G_q), \qquad KK^{G_q}_1(C(G_q/T), \mathbb{C}) = 0.
$$
\end{cor}
Let us also discuss the following result which is closely related to theorem \ref{PD}.
\begin{theorem} \label{dsdouble}
The standard Podle\'s sphere $ C(G_q/T) $ is a direct summand of $ \mathbb{C} \oplus \mathbb{C} $ in $ KK^{\DD(G_q)} $.
\end{theorem}
\proof Let us consider the elements $ \alpha_q \in KK^{\DD(G_q)}_0(C(G_q/T), \mathbb{C} \oplus \mathbb{C}) $ and $ \beta_q
\in KK^{\DD(G)}_0(\mathbb{C} \oplus \mathbb{C}, C(G_q/T)) $ given by
$$
\alpha_q = [D] \oplus [D \otimes \E_{-1}], \qquad \beta_q = (-[\E_1]) \oplus [\E_0],
$$
respectively. Following the proof of theorem \ref{PD} we shall show $ \alpha_q \circ \beta_q = \id $. \\
Consider first the case $ q = 1 $. All $ C^*(G_1) $-coactions in the construction of $ \alpha_1 $ and $ \beta_1 $ are trivial,
and it suffices to check $ \alpha_1 \circ \beta_1 = \id $ in $ KK^{G_1}_0(C(G_1/T), C(G_1/T)) $.
Using $ KK^{G_1}_*(C(G_1/T), C(G_1/T)) \cong KK^T_*(C(G_1/T), \mathbb{C}) $
and the universal coefficient theorem for $ KK^T $ we only have
to compare the corresponding actions on $ K^T_*(C(G_1/T)) $.
One obtains
\begin{align*}
[\E_0] \circ \alpha_1 \circ \beta_1 \circ \epsilon &= [\E_0] \circ [D \otimes \E_{-1}] \circ [\E_0] \circ \epsilon - [\E_0] \circ [D]
\circ [\E_1] \circ \epsilon \\
&= [\E_0] \circ [D \otimes \E_{-1}] = [\mathbb{C}_0]
\end{align*}
in $ R(T) $ due to proposition \ref{indexhomogenous}, and similarly $ [\E_1] \circ \alpha_1 \circ \beta_1 \circ \epsilon = [\mathbb{C}_1] $.
Taking into account McLeod's theorem \cite{McLeod} this yields the assertion for $ q = 1 $. \\
For general $ q \in (0,1] $ we recall
$$
KK^{\DD(G_q)}_*(C(G_q/T), C(G_q/T)) \cong KK^T_*(C(G_q/T) \twisted_{G_q} C(G_q), \mathbb{C})
$$
and notice that the elements corresponding to $ \alpha_t \circ \beta_t $ for $ t \in [q,1] $ assemble into a class
in $ KK^T(C({\bf G}/T) \twisted_{{\bf G}} C({\bf G}), C[q,1]) $. The comparison argument in the proof of
theorem \ref{PD} carries over and yields $ \alpha_q \circ \beta_q = \id $ in $ KK^{\DD(G_q)} $. \qed \\
On the level of $ G_q $-equivariant Kasparov theory one can strengthen
the assertion of theorem \ref{dsdouble} as follows.
\begin{prop} \label{podleskkg}
The standard Podle\'s sphere $ C(G_q/T) $ is isomorphic to $ \mathbb{C} \oplus \mathbb{C} $ in $ KK^{G_q} $.
\end{prop}
\proof We have already seen that the elements $ \alpha_q $ and $ \beta_q $ defined in theorem \ref{dsdouble} satisfy $ \alpha_q \circ \beta_q = \id $
in $ KK^{\DD(G_q)} $, hence this relation holds in $ KK^{G_q} $ as well.
Using proposition \ref{indexhomogenous} one immediately calculates $ \beta_q \circ \alpha_q = \id $ in $ KK^{G_q} $. \qed

\bibliographystyle{plain}

\end{document}